\documentclass[a4paper, 12pt, leqno]{article}

\usepackage{amsmath,amsfonts}

\usepackage[active]{srcltx}
\begin{document}
\bibliographystyle{plain}
\newtheorem{theo}{Theorem}[section]
\newtheorem{lemme}[theo]{Lemma}
\newtheorem{cor}[theo]{Corollary}
\newtheorem{defi}[theo]{Definition}
\newtheorem{prop}[theo]{Proposition}
\newtheorem{problem}[theo]{Problem}
\newtheorem{remarque}[theo]{Remark}
\newcommand{\beq}{\begin{eqnarray}}
\newcommand{\enq}{\end{eqnarray}}
\newcommand{\be}{\begin{eqnarray*}}
\newcommand{\en}{\end{eqnarray*}}
\newcommand{\Td}{\mathbb T^d}
\newcommand{\Rd}{\mathbb R^n}
\newcommand{\R}{\mathbb R}
\newcommand{\N}{\mathbb N}
\newcommand{\Sn}{\mathbb S}
\newcommand{\Zd}{\mathbb Z^d}
\newcommand{\Linf}{L^{\infty}}
\newcommand{\dt}{\partial_t}
\newcommand{\Dt}{\frac{d}{dt}}
\newcommand{\Dtt}{\frac{d^2}{dt^2}}
\newcommand{\demi}{\frac{1}{2}}
\newcommand{\vf}{\varphi}
\newcommand{\epu}{_{\epsilon}}
\newcommand{\ep}{^{\epsilon}}
\newcommand{\bfi}{{\mathbf \Phi}}
\newcommand{\bpsi}{{\mathbf \Psi}}
\newcommand{\bx}{{\mathbf x}}
\newcommand{\ds}{\displaystyle}
\newcommand{\mT}{\mathfrak{T}}

\let\cal=\mathcal
\title{On the regularity of solutions \\ of optimal transportation problems}
\author{Gr\'egoire Loeper\footnotemark[1]}
\maketitle
\footnotetext[1]{Institut Camille Jordan, Universit\'e Claude Bernard Lyon 1}

\begin{abstract}
We give a necessary and sufficient condition on the cost function so that the map solution of Monge's optimal  transportation problem is continuous for arbitrary smooth positive data.
This condition was first introduced by Ma, Trudinger and Wang \cite{MTW, TW} for a priori estimates of the corresponding Monge-Amp\`ere equation. It is expressed by a so-called {\em cost-sectional curvature} being non-negative. We show that when the cost function is the squared distance of  a Riemannian manifold, the cost-sectional curvature yields the sectional curvature. As a consequence,  if the manifold does not have non-negative sectional curvature everywhere, the optimal transport map {\em cannot be continuous} for arbitrary smooth positive data.
The non-negativity of the cost-sectional curvature is shown to be equivalent to the connectedness of the contact set between any cost-convex function (the proper generalization of a convex function) and any of its supporting functions.
When the cost-sectional curvature is uniformly positive,  we obtain that optimal maps are continuous or H\"older continuous under quite weak assumptions on the data, compared to what is needed in the Euclidean case. This case includes the quadratic cost on the round sphere.

\end{abstract}

\newpage
\tableofcontents

\section{Introduction}

Given $A,B$ two topological spaces, a cost function  $c:A \times B \to \R$,  and $\mu_0, \mu_1$ two probability measures respectively on $A$ and $B$,  Monge's problem of optimal transportation consists in finding among all measurable maps $T:A\to B$  that push forward $\mu_0$ onto $\mu_1$ (hereafter $T_{\#}\mu_0 = \mu_1$) 
in the sense that 
\beq\label{pushforward}
\forall E \subset B \text{ Borel }, \mu_1(E)=\mu_0(T^{-1}(E)),
\enq
a map that realizes
\beq\label{monge}
{\rm Argmin} \Big\{\int_A c(x,T(x))\,d\mu_0(x), T_{\#}\mu_0=\mu_1\Big\}.
\enq 
Optimal transportation has undergone a rapid and important development since the pioneering work of Brenier, who discovered that when $A=B=\Rd$ and the cost is the distance squared, optimal maps for the problem (\ref{monge}) are gradients of convex functions \cite{Br1} (see also \cite{Knott} where the connection with gradients was first proved). Following this result and its subsequent extensions, the theory of optimal transportation has flourished, with generalizations to other cost functions \cite{CaGenCost, McGa}, more general spaces such as Riemannian manifolds \cite{Mc2},  applications in many other areas of mathematics such as  geometric analysis, functional inequalities, fluid mechanics, dynamical systems, and other more concrete applications such as irrigation, cosmology. 

When $A,B$ are domains of the Euclidean space $\Rd$, or of a Riemannian manifold,  a  common feature to all optimal transportation problems is that optimal maps derive from a (cost-convex) potential, which, assuming some smoothness,  is in turn solution to a fully non-linear elliptic PDE: the Monge-Amp\`ere equation. In all cases, the Monge-Amp\`ere equation arising from an optimal transportation problem reads in local coordinates
\beq\label{magegen}
\det (D^2\phi-{\cal A}(x,\nabla\phi))=f(x,\nabla\phi),
\enq
where $(x,p)\to{\cal A}(x,p)$ is a symmetric matrix valued function, that depends on the cost function $c(x,y)$ through the formula 
\beq\label{defA}
{\cal A}(x,p)=-D^2_{xx}c(x,y) \text{ for } y \text{ such that }-\nabla_xc(x,y)=p.
\enq
That there is indeed a unique $y$ such that $-\nabla_xc(x,y)=p$ will be guaranteed by condition {\bf A1} given hereafter. The optimal map will then be $$x\to y: -\nabla_xc(x,y)=\nabla\phi(x).$$ 
 In the case ${\cal A}=0$,  equation (\ref{magegen}) was well known and studied before optimal transportation since it appears in Minkowsky's problem: find a convex hypersurface with prescribed Gauss curvature. In the case of optimal transportation, the boundary condition consists in prescribing that the image of the optimal map equals a certain domain. It is known as the second boundary value problem.

Until recently, except in the  particular case of the so-called reflector antenna, treated by Wang \cite{Wa4} (see also \cite{CaGuHu} for $C^1$ regularity),  the regularity of optimal maps was only known in the case where the cost function is the (Euclidean) squared distance  $c(x,y)=|x-y|^2$, which is the cost considered by Brenier in \cite{Br1}, for which the matrix ${\cal A}$ in (\ref{magegen}) is the identity (which is trivially equivalent to the case  ${\cal A}=0$). Those results have involved several authors, among which Caffarelli, Urbas , and Delano\"e. An important step was made recently by Ma, Trudinger and Wang  \cite{MTW}, and Trudinger and Wang  \cite{TW}, who introduced a condition (named A3 and   A3w in their papers) on the cost function under which they could show existence of smooth solutions to (\ref{magegen}). Let us give right away this condition that will play a central role in the present paper.
Let $A=\Omega, B=\Omega'$ be bounded domains of $\Rd$ on which the initial and final measures will be supported. Assume that $c$ belongs to   $C^4(\Omega\times\Omega')$. For $(x,y) \in (\Omega \times \Omega'), (\xi, \nu) \in \Rd\times\Rd$,  we define
\beq\label{R}
&&\mathfrak{S_c}(x,y)(\xi,\nu):= D^2_{p_kp_l}{\cal A}_{ij}\,\xi_i\xi_j\,\nu_k\nu_l\,(x,p),  \  \  p=-\nabla_xc(x,y).
\enq
Whenever $\xi,\nu$ are orthogonal unit vectors,  we will say that $\mathfrak{S_c}(x,y)(\xi,\nu)$ defines the {\em cost-sectional curvature from $x$ to $y$ in the directions $(\xi, \nu)$}.
As we will see in Definition \ref{defiR}, this definition is intrinsic.
Note that this map is in general not symmetric, and that it depends on two points $x$ and $y$.
The reason why we use the word sectional curvature will be clear in a few lines. 
We will say that the cost function  $c$ has non-negative cost-sectional curvature on $(\Omega\times\Omega')$, if
\beq\label{A}
\mathfrak{S_c}(x,y)(\xi,\nu)\geq 0  \  \ \forall (x,y) \in (\Omega \times \Omega'), \forall (\xi, \nu) \in \Rd\times\Rd, \  \xi\perp\nu.
\enq
A cost function satisfies condition {\bf Aw} on $(\Omega\times\Omega')$ if and only if it has non-negative cost-sectional curvature on $(\Omega\times\Omega')$, i.e. if it satisfies (\ref{A}).

Under condition {\bf Aw} and natural requirements on the domains $\Omega,\Omega'$,  Trudinger and Wang \cite{TW} showed  that the solution to (\ref{magegen}) is globally smooth for smooth positive measures $\mu_0, \mu_1$. They showed that  {\bf Aw} is satisfied by a large class of cost functions, that we will give as examples later on. Note that the quadratic cost satisfies  assumption {\bf Aw}. This result is achieved by the so-called continuity method, for which a key ingredient is to obtain a priori estimates on the second derivatives of the solution. At this stage, condition {\bf Aw} was used in a crucial way. 
However, even if it was known that not all cost functions can lead to smooth optimal maps, it was unclear whether the condition {\bf Aw} was necessary, or just a technical condition for the a-priori estimates to go through.

In this paper we show that the condition {\bf Aw} is indeed the {\em necessary and sufficient condition for regularity:} one can not expect regularity without this condition, and more precisely, if $\mathfrak{S_c}(x,y)(\xi,\nu)<0$ for $(x,y)\in (\Omega\times\Omega'), \xi\perp\nu \in \Rd$, one can immediately build a pair of $C^\infty$ strictly positive measures, supported on sets that satisfy the usual smoothness and convexity assumptions, so that the optimal potential is not even $C^1$, and the optimal map is therefore discontinuous.
This result is obtained by analyzing the geometric nature of condition (\ref{A}). Let us first recall that the solution $\phi$ of the Monge-Amp\`ere equation is a priori known to be cost-convex (in short c-convex), meaning that at each point $x\in \Omega$, there exist $y\in \Omega'$ and a value $\phi^c(y)$  such that 
\be 
&&-\phi^c(y)-c(x,y)=\phi(x),\\
&&-\phi^c(y)-c(x',y) \leq \phi(x'), \ \forall x'\in \Omega.
\en
The function $-\phi^c(y) - c(x,y)$ is called a supporting function, and the function $y\to \phi^c(y)$ is called the cost-transform (in short the c-transform) of $\phi$, also defined by
\be
\phi^c(y)=\sup_{x\in \Omega}\{-c(x,y)-\phi(x)\}.
\en
(These notions will be recalled in greater details hereafter.)
We prove that the condition {\bf  Aw} can be reformulated as a property of cost-convex functions, which we call {\em connectedness of the contact set}:
\beq\label{C}
&&\text{For all } x\in \Omega,\text{ the contact set }G_\phi(x):= \{y: \phi^c(y)=-\phi(x)-c(x,y)\}\\
&& \text{ is connected. } \nonumber
\enq
Assuming a natural condition on $\Omega'$  (namely its c-convexity, see Definition \ref{def-c-convex}) this condition involves only the cost function since it must hold for any  $\phi^c$ defined through a c-transform.

A case of special interest for applications is the generalization of Brenier's cost $\demi|x-y|^2$ to  Riemannian manifolds, namely $c(x,y)=\demi d^2(x,y)$. Existence and uniqueness of optimal maps in that case was established by McCann \cite{Mc2}, and further examined by several authors, with many interesting applications in geometric and functional analysis (for example \cite{cormac, OV}). The optimal map  takes the form $x\to \exp_x(\nabla\phi(x))$ for $\phi$ a c-convex potential and is called a gradient map. Then, a natural question is the interpretation of condition {\bf Aw} and of the cost-sectional curvature in this context. We show that  for some universal constant $K$, 
\be
\text{Cost-sectional curvature from } x \text{ to } x= K\cdot\text{ Riemannian sectional curvature at } x.
\en
(We mean there that the equality holds for every 2-plane and actually $K=2/3$.)
As a direct consequence of the previous result,  {\em the optimal (gradient) map will  not be continuous for arbitrary smooth positive data if the manifold does not have non-negative sectional curvature everywhere.}
Although the techniques are totally different, it is interesting to notice that in recent works, Lott \& Villani \cite{LV}, and Sturm \cite{St} have recovered the Ricci curvature through a property of optimal transport maps (namely through the displacement convexity of some functionals). Here, we somehow recover the sectional curvature through the continuity of optimal maps.

We next investigate the continuity of optimal maps under the stronger condition of uniformly positive cost-sectional curvature, or condition {\bf As}:
\beq\label{AA}
\exists C_0>0: \  \mathfrak{S_c}(x,y,\xi,\nu) \geq C_0|\xi|^2|\nu|^2,  \  \  \   \forall (x,y) \in (\Omega \times \Omega'), (\xi, \nu) \in \Rd\times\Rd,  \xi\perp\nu.
\enq
We obtain that the  (weak) solution of (\ref{magegen}) is $C^1$ or $C^{1,\alpha}$ under quite mild assumptions on the measures. Namely, for $B_r(x)$  the ball of radius $r$ and center $x$,  $\mu_1$ being bounded away from 0, we need $\mu_0(B_r(x))=o(r^{n-1})$
to show that the solution of (\ref{magegen}) is $C^1$ and $\mu_0(B_r(x))=O(r^{n-p}), p<1$ to show that it is $C^{1,\alpha}$, for $\alpha=\alpha(n,p)\in (0,1)$.
Those conditions allow $\mu_0,\mu_1$ to be singular with respect to the Lebesgue measure and $\mu_0$  to vanish. 

This result can be seen as analogous to Caffarelli's $C^{1,\alpha}$ estimate \cite{Ca2} for a large class of  cost functions and related Monge-Amp\`ere equations. It also shows that the partial regularity results are better under {\bf As} than under {\bf Aw}, since Caffarelli's $C^{1,\alpha}$ regularity result required $\mu_0,\mu_1$ to have densities bounded away from 0 and infinity, and it is  known to be close to optimal \cite{Wa2}.

In a forthcoming work \cite{Lsphere} we shall prove that the quadratic cost on the sphere has uniformly positive cost-sectional curvature, i.e. satisfies {\bf As}. We obtain therefore regularity of optimal (gradient) maps under adequate conditions.

The rest of the paper is organized as follows:
in section 2 we gather all definitions and results that we will need throughout the paper. In section 3 we state our results.
Then each following section is devoted to the proof of a theorem. 
The reader knowledgeable about the subject might skip directly to section 3.

\section{Preliminaries}

\subsection{Notation}
Hereafter ${\rm dVol}$ denotes the Lebesgue measure of $\Rd$ 
and $B_r(x)$ denotes a ball of radius $r$ centered at $x$. For $\delta>0$, we set classically $\Omega_\delta=\{x\in \Omega, d(x,\partial\Omega)>\delta)\}$.  When we say that a function (resp. a measure) is smooth without stating the degree of smoothness, we assume that it is $C^\infty$-smooth  (resp. has a $C^\infty$-smooth density with respect to the Lebesgue measure).

\subsection{Kantorovitch duality and c-convex potentials}\label{section-c-convex}
In this section, we recall how to obtain the optimal map from a c-convex potential in the general case. This allows us to introduce definitions that we will be using throughout the paper. References concerning the existence of optimal map by Monge-Kantorovitch duality are \cite{Br1} for the cost $|x-y|^2$, \cite{McGa} and \cite{CaGenCost} for general costs, \cite{Mc2} for the Riemannian case, otherwise the book \cite{Vi} offers a rather complete reference on the topic.

Monge's problem (\ref{monge}) is first relaxed to become a problem of linear programming; one seeks now 
\beq\label{kanto*}
{\cal I}=\inf\Big\{\int_{\Rd\times\Rd} c(x,y)d\pi(x,y); \  \  \pi \in \Pi(\mu_0,\mu_1)\Big\}
\enq
where $\Pi(\mu_0,\mu_1)$ is the set of positive measures on $\Rd\times\Rd$ whose marginals are respectively $\mu_0$ and $\mu_1$. Note that the (Kantorovitch) infimum (\ref{kanto*}) is smaller than the (Monge) infimum of the cost (\ref{monge}), since whenever a map $T$ pushes forward $\mu_0$ onto $\mu_1$, the measure $\pi_T(x):= \mu_0(x)\otimes\delta_{T(x)}(y)$ belongs to $\Pi(\mu_1,\mu_1)$.

Then, the dual Monge-Kantorovitch problem is to find an optimal pair of potentials $(\phi, \psi)$ that realizes 
\beq\label{kanto}
 \hspace{1 cm} {\cal J} = \sup\Big\{ -\int \phi(x)d\mu_0(x) - \int \psi(y)d\mu_1(y); \phi(x)+\psi(y) \geq -c(x,y)\Big\} .
\enq

The constraint on $\phi, \psi$ leads to the definition of  c(c*)-transforms:

\begin{defi} Given a lower semi-continuous function $\phi: \Omega\subset \Rd\to \R\cup\{+\infty\}$, we define its c-transform at  $y\in\Omega'$ by
\beq\label{def-c-transform}
\phi^c(y)=\sup_{x\in \Omega}\{- c(x,y) - \phi(x)\}.
\enq
Respectively, for $\psi:\Omega'\subset\Rd \to \R$ also lower semi-continuous, define its c*-transform at $x\in\Omega$ by
\beq\label{def-c*-transform}
\psi^{c*}(x)=\sup_{y\in \Omega'}\{- c(x,y) - \psi(y)\}.
\enq
A function is said cost-convex, or, in short, c-convex, if it is the c*-transform of another function $\psi:\Omega' \to \R$, i.e. for $x \in \Omega$, $\phi(x) = \sup_{y\in\Omega'} \{-c(x,y) - \psi(y)\}$, for some lower semi-continuous $\psi: \Omega' \to \R$.  Moreover in this case $\phi^{cc*}:=(\phi^c)^{c*}=\phi$ on $\Omega$ (see \cite{Vi}).
\end{defi}

Our first assumption on $c$ will be:
\begin{enumerate}
\item[{\bf A0}] The cost-function $c$ belongs to  $C^4(\bar\Omega\times\bar\Omega')$. \end{enumerate}
We will also always assume  that $\Omega, \Omega'$ are bounded. 
These assumptions are not the weakest possible for the existence/uniqueness theory.
\begin{prop}\label{prop-c-convex}
If $c$ is Lipschitz and semi-concave with respect to $x$, locally uniformly with respect to $y$, and
if $\Omega'$ is bounded, $\phi^c$  will be locally semi-convex and Lipschitz.
In  particular, this holds under assumption {\bf A0}. The symmetric statement holds for $\psi^{c*}$.  
\end{prop}
By Fenchel-Rockafellar's duality theorem, we have ${\cal I}={\cal J}$.  One can then easily show that the supremum (\ref{kanto})  and the infimum (\ref{kanto*}) are achieved. 
Since the condition $\phi(x)+\psi(y) \geq -c(x,y)$ implies $\psi \geq \phi^c$, we can assume that for the optimal pair in ${\cal J}$ we have $\psi=\phi^c$ and $\phi=\phi^{cc*}$. Writing the equality of the integrals in (\ref{kanto*}, \ref{kanto}) for any optimal $\gamma$ and any optimal pair $(\phi, \phi^c)$ we obtain that $\gamma$ is supported in
$\Big\{\phi(x)+\phi^c(y)+c(x,y)=0\Big\}$. 
This leads us to the following definition:
\begin{defi}[Gradient mapping]\label{def-grad-map}
Let $\phi$ be a c-convex function. We define the set-valued mapping $G_\phi$ by 
\be
G_\phi(x)=\Big\{y\in \Omega', \phi(x)+\phi^c(y)=-c(x,y)\Big\}.
\en
For all $x\in \Omega$, $G_\phi(x)$ is the contact set between $\phi^c$ and its supporting function $-\phi(x)-c(x,\cdot)$.  
\end{defi}
Noticing that for all $y\in G_\phi(x)$, $\phi(\cdot) + c(\cdot, y)$ has a global minimum at $x$, we introduce /  recall  the following definitions:
\begin{defi}[subdifferential]\label{sub-dif}
For $\phi$ a semi-convex function, the subdifferential of $\phi$ at $x$, that we denote $\partial\phi(x)$, is the set 
\be
\partial\phi(x)=\Big\{p\in \Rd, \phi(y) \geq \phi(x) + p\cdot(y-x) + o(|x-y|)\Big\}.
\en
\end{defi}
The subdifferential is always a convex set, and is always non empty for a semi-convex function. 

\begin{defi}[c-subdifferential]\label{c-sub-dif}
If $\phi$ is $c$-convex, the c-sub-differential of $\phi$ at $x$, that we denote $\partial^c\phi(x)$, is the set
\be
\partial^c\phi(x)=\Big\{ -\nabla_x c(x,y), y\in G_\phi(x)\Big\}.
\en
The inclusion $\emptyset \neq \partial^c\phi(x)\subset \partial\phi(x)$ always holds.
\end{defi}
We introduce now two  assumptions on the cost-function, which are the usual assumptions made in order to obtain an optimal map. For $x=(x_1,...,x_n), y=(y_1...y_n)$, let us first introduce the notation 
\be 
D^2_{xy}c(x,y) = \left[\partial_{x_i}\partial_{y_j} c(x,y)\right]_{1 \leq i,j \leq n}.
\en

\begin{enumerate}
\item[{\bf A1}] For all $x\in \bar\Omega$, the mapping $y\to -\nabla_xc(x,y)$ is injective on $\bar\Omega'$. 
\item[{\bf A2}] The cost function $c$ satisfies $\det D^2_{xy}c \neq 0$ for all $(x,y)\in \bar\Omega\times\bar\Omega'$.
\end{enumerate}
This leads us to the definition of the {\em c-exponential map}:
\begin{defi}\label{defiT} 
Under assumption {\bf A1}, for $x\in \Omega$ we define the c-exponential map at $x$, which we denote by $\mathfrak{T}_x$, such that 
\be
\forall (x,y)\in (\Omega\times\Omega'), \mathfrak{T}_x(-\nabla_xc(x,y))=y.
\en
Moreover, under assumptions {\bf A0}, {\bf A1}, {\bf A2}, and assuming that $\Omega'$ is connected, there exists a constant $C_\mT>0$ 
that depends on $c, \Omega, \Omega'$,  such that for all $x \in \Omega$, for all $p_1, p_2 \in -\nabla_x c(x, \Omega')$, 
\beq\label{bma}
\frac{1}{C_\mT} \leq \frac{|\mathfrak{T}_{x}(p_2) - \mathfrak{T}_{x}(p_1)|}{|p_2-p_1|} \leq C_\mT.
\enq

\end{defi}

\textsc{Remark 1.} The definition c-exponential map is again motivated by the case cost=distance squared, where the c-exponential map is the exponential map.
Moreover, notice the important identity
\beq\label{dpexp}
[D^2_{xy}c]^{-1}=-D_p\mathfrak{T}_x\big|_{x, p=-\nabla_xc(x,y)}.
\enq

\textsc{Remark 2.} Anticipating the extension to Riemannian manifolds, we mention at this point that this definition is intrinsic, i.e. it defines in a coordinate independent way the map $\mathfrak{T}$ as a map going from $M\times TM$ to $M$. In this setting, the gradients should be computed with respect to the metric $g$ of the manifold.

Under assumptions {\bf A1}, {\bf A2}, $G_\phi$ is single valued outside of a set of Hausdorff dimension less than or equal to $n-1$, hence, 
if $\mu_0$ does not give mass to sets of Hausdorff dimension less than $n-1$, $G_\phi$ will be the optimal map for Monge's problem while the optimal measure in (\ref{kanto*}) will  be $\pi = \mu_0\otimes\delta_{G_\phi(x)}$. So, after having relaxed the constraint that the optimal $\pi$ should be  supported on the graph of a map, one still obtains a minimizer that satisfy this  constraint.

Notice  that Monge's historical cost was equal to the distance itself: $c(x,y)=|x-y|$. One sees immediately that for this cost function, there is not a unique $y$ such that $-\nabla_xc(x,y)=\nabla\phi(x)$, hence assumption {\bf A1} is not satisfied and, indeed, there is in general no uniqueness of the optimal map.

We now state a general existence theorem, under  assumptions that are clearly not minimal, but that will suffice for the scope of this paper, where we deal with regularity issues.
\begin{theo}
Let $\Omega, \Omega'$ be two bounded domains of $\Rd$. Let $c\in C^4(\bar\Omega\times\bar\Omega')$ satisfy assumptions {\bf A0}-{\bf A2}. Let $\mu_0, \mu_1$ be two probability measures on $\Omega$ and $\Omega'$. Assume that $\mu_0$ does not give mass to sets of Hausdorff dimension less than or equal to $n-1$. Then there exists a $d\mu_0$ a.e. unique minimizer $T$ of Monge's optimal transportation problem (\ref{monge}). Moreover, there exists $\phi$ c-convex on $\Omega$ such that $T=G_\phi$  (see \ref{def-grad-map}). Finally, if $\psi$ is c-convex and satisfies $G_{\psi\,\#}\mu_0=\mu_1$, then $\nabla\psi=\nabla\phi$ $d\mu_0$ a.e.
\end{theo}

\subsection{Notion of c-convexity for sets}
Following \cite{MTW, TW}, we introduce here the notions that extend naturally the notions of  convexity / strict convexity for a set. 

\begin{defi}[c-segment]\label{def-c-seg}
Let $p \to \mathfrak{T}_x(p)$ be the mapping defined by assumption {\bf A1}. The point $x$ being held fixed, a c-segment with respect to $x$ is the image by $\mathfrak{T}_x$ of a segment of $\Rd$. 

If for $v_0, v_1 \in \Rd$ we have $\mathfrak{T}_x(v_i)=y_i, i=0,1$, the c-segment with respect to $x$ joining $y_0$ to $y_1$ will be $\{y_\theta, \theta\in [0,1]\}$ where $y_\theta=\mathfrak{T}_x(\theta v_1+ (1-\theta) v_0)$. It will be denoted $[y_0,y_1]_x$.
\end{defi}

\begin{defi}[c-convex sets]\label{def-c-convex}
Let $\Omega,\Omega' \subset \Rd$. We say that $\Omega'$ is $c$-convex with respect to $\Omega$ if for all $y_0,y_1 \in \Omega', x\in \Omega$, the c-segment $[y_0, y_1]_x$ is contained in $\Omega'$.
\end{defi}

\textsc{Remark.} Note that this can be said in the following way:
for all $x\in \Omega$, the set $-\nabla_xc(x,\Omega')$ is convex.

\begin{defi}[uniform strict c-convexity of sets]\label{def-unif-c-convex}
For $\Omega, \Omega'$ two subsets of $\Rd$, we say that $\Omega'$ is uniformly strictly c-convex with respect to $\Omega$ if the sets $\{-\nabla_xc(x,\Omega')\}_{x\in \Omega}$ are uniformly strictly convex, uniformly with respect to $x$. We say that $\Omega$ is uniformly strictly c*-convex with respect to $\Omega'$ if the dual assertion holds true.
\end{defi}

\textsc{Remark 1.} In local coordinates, $\Omega$ is uniformly strictly c*-convex with respect to $\Omega'$  reads
\beq\label{gamma}
[D_i\gamma_j(x) - D_{p_k}{\cal A}_{ij}(x,p)\gamma_k]\tau_i\tau_j \geq \epsilon_0 >0,
\enq
for some $\epsilon_0>0$, for all $x\in \partial\Omega, p\in  -\nabla_xc(x,\Omega')$, unit tangent  vector  $\tau$ and outer unit normal $\gamma$.

\textsc{Remark 2.} When ${\cal A}$ does not depend on $p$, one recovers the usual convexity.

\paragraph{Remarks on the sub-differential and c-sub-differential}
The question is to know if we have for all $\phi$ c-convex on $\Omega$, for all $x\in \Omega$, $\partial\phi(x)=\partial^c\phi(x)$. 
Clearly, when $\phi$ is c-convex and differentiable at $x$,the equality holds.
For $p$ an extremal point of $\partial\phi(x)$, there will be a sequence $x_n$ converging to $x$ such that $\phi$ is differentiable at $x_n$ and $\lim_n\nabla\phi(x_n)=p$. Hence, extremal points of $\partial\phi(x)$ belong to $\partial^c\phi(x)$. Then it is not hard to show the
\begin{prop}\label{remarks}
Assume that $\Omega'$ is c-convex with respect to $\Omega$. The following assertions are equivalent:
\begin{enumerate}
\item For all $\phi$ c-convex on $\Omega$, $x\in \Omega$,  $\partial^c\phi(x)=\partial\phi(x)$.
\item  For all $\phi$ c-convex on $\Omega$, $x\in \Omega$,  $\partial^c\phi(x)$ is convex.
\item For all $\phi$ c-convex on $\Omega$, $x\in \Omega$, $G_\phi(x)$ is c-convex with respect to $x$.
\item For all $\phi$ c-convex on $\Omega$, $x\in \Omega$, $G_\phi(x)$ is connected.
\end{enumerate}
\end{prop}

\textsc{Proof.} We prove only that (4) implies (2). First, the connectedness of $G_\phi(x)$ implies the connectedness of $\partial^c\phi(x)$, since $\nabla_xc$ is continuous. Then for $x_0\in\Omega$, $y_0, y_1 \in \Omega'$, assume that $y_0$ and $y_1$ both belong to $G_\phi(x_0)$. Letting $$h(x)= \max\{-c(x,y_0) + c(x_0, y_0) + \phi(x_0), -c(x,y_1) + c(x_0, y_1) + \phi(x_0)\},$$ one has $\phi(x)\geq h(x)$ on $\Omega$, with equality at $x=x_0$. 
Hence $\partial^c h(x_0)\subset \partial^c \phi(x_0)$. Since the property (4) is satisfied, $\partial^c h(x_0)$ is connected, and as it is included in $\partial h(x_0)$ which is a segment, it is equal to the segment $[-\nabla_xc(x_0, y_0), -\nabla_xc(x_0, y_1)]$. This shows that $\partial^c\phi(x_0)$ is convex. 

$\hfill \Box$ 

\subsection{The Monge-Amp\`ere equation}
In all cases, for $\phi$ a $C^2$ smooth c-convex potential  such that
$G_{\phi\,\#} \mu_0  = \mu_1$,
the conservation of mass is expressed in local coordinates by the following Monge-Amp\`ere equation 
\beq\label{ma0} 
\det(D^2_{xx}c(x,G_\phi(x)) + D^2\phi) = |\det D^2_{xy}c|\frac{\rho_0(x)}{\rho_1(G_\phi(x))},
\enq
where  $\rho_i={\rm d}\mu_i/{\rm dVol}$ denotes the density of $\mu_i$ with respect to the Lebesgue measure. (See \cite{MTW} for a derivation of this equation, or \cite{cormac}, \cite{del2}.) Hence, the equation fits into the general form (\ref{magegen}).

\subsection{Generalized solutions}
\begin{defi}[Generalized solutions]\label{defiweak}
Let $\phi:\Omega\to \R$ be a c-convex function. Then
\begin{itemize}
\item $\phi$ is a weak Alexandrov solution to (\ref{ma0}) if and only if
\beq\label{weakalex}
\text{for all } B\subset\Omega, \ \mu_0(B)=\mu_1(G_\phi(B)).
\enq
This will be denoted by $\mu_0 = G_\phi^{\#}\mu_1$.

\item  $\phi$ is a weak Brenier solution to (\ref{ma0}) if and only if 
\beq\label{weakbr}
\text{for all } B'\subset\Omega', \ \mu_1(B')=\mu_0(G_\phi^{-1}(B')).
\enq
This is equivalent to $\mu_1 = G_{\phi\,\#}\mu_0$.

\end{itemize}

\end{defi}

\paragraph{Alexandrov and Brenier solutions}

First notice that in the definition (\ref{weakbr}), $\mu_1$ is deduced from $\mu_0$, while it is  the contrary in (\ref{weakalex}).
As we have seen, the Kantorovitch procedure (\ref{kanto}) yields an optimal transport map whenever $\mu_0$ does not give mass to sets of Hausdorff dimension less than $n-1$. Moreover, the map $G_\phi$ will satisfy (\ref{weakbr}) by construction, and hence will be a weak Brenier solution to (\ref{ma0}).  Taking advantage of the c-convexity of $\phi$ one can show that whenever $\mu_1$ is absolutely continuous with respect to the Lebesgue measure, $G_\phi^{\#}\mu_1$ is countably additive, and hence is a Radon measure (see \cite[Lemma 3.4]{MTW}); then a Brenier solution is an Alexandrov solution.
Note that one can consider $\mu_0=G_\phi^{\#}{\rm dVol}$, this will be the Monge-Amp\`ere measure of $\phi$. Most importantly, for $\mu_0$ supported in $\Omega$, $G_{\phi\,\#}\mu_0 = {\bf 1}_{\Omega'}{\rm dVol}$ does not imply that $G_\phi^{\#}{\rm dVol}=\mu_0$, except if $\Omega'$ is c-convex with respect to $\Omega$  (see \cite{MTW}).

\subsection{Cost-sectional curvature and conditions {\bf Aw}, {\bf As}}
\label{conditions}

A central notion in the present paper will be the notion of {\em cost-sectional curvature} $\mathfrak{S_c}(x,y)$. 
\begin{defi}\label{defiR}
Under assumptions {\bf A0}-{\bf A1}-{\bf A2},  one can define  on $T_x\Omega \times T_x\Omega$ the real-valued map
\beq\label{defR}
 \  \  \  \  \   \mathfrak{S_c}(x_0,y_0)(\xi,\nu)= D^4_{p_\nu p_\nu x_\xi x_\xi}\Big[(x,p)\to -c(x,\mathfrak{T}_{x_0}(p))\Big]\Big|_{x_0,p_0=-\nabla_x c(x_0,y_0)}.
\enq
When $\xi,\nu$ are unit orthogonal  vectors, $\mathfrak{S_c}(x_0,y_0)(\xi,\nu)$ defines  the  cost-sectional curvature from $x_0$ to $y_0$ in directions $(\xi,\nu)$.
The definition (\ref{defR}) is equivalent to the following:
\beq\label{defRint}
 \  \  \  \  \   \mathfrak{S_c}(x_0,y_0)(\xi,\nu)= D^2_{tt}D^2_{ss} \Big[(s,t)\to -c(\exp_{x_0}(t\xi), \mathfrak{T}_{x_0}(p_0+s\nu))\Big]\Big|_{t,s=0}.
\enq
\end{defi}
The fact that the definition (\ref{defRint}) and (\ref{defR}) are equivalent follows from the following observation:

\begin{prop}\label{propnewR}
The definition of $\mathfrak{S_c}(x_0,y_0)(\xi,\nu)$ is intrinsic, i.e. depends only on $(x_0, y_0)\in \Omega\times\Omega'$ and on  $(\xi, \nu)\in T_{x_0}(\Omega) \times T_{x_0}(\Omega)$, and not on the choice of local coordinates around $x_0$ or $y_0$.
Moreover, it is symmetric: letting $c^*(y,x) = c(x,y)$, and $\mathfrak{T}^*$ be the c*-exponential map, the identity
\beq
\label{newRint}\mathfrak{S_c}(x_0,y_0)(\xi,\nu) = \mathfrak{S_{c*}}(y_0,x_0)(\tilde\nu,\tilde\xi)
\enq
holds with $\tilde\nu = D_{p} \mathfrak{T}_{x_0}(p_0)\cdot \nu$,and $\tilde\xi = [D_{q}\mathfrak{T}^*_{y_0}(q_0)]^{-1}\cdot\xi$, with $p_0$ as above and 
$q_0=-\nabla_y c(x_0,y_0)$.
Notice that whenever $\xi\perp\nu$, one has $\tilde\xi\perp\tilde\nu$.
\end{prop}
 
\textsc{Proof.} The proof is deferred to the appendix.

\textsc{Remark.} The intrinsic nature of the cost-sectional curvature tensor has been observed independently in \cite{KimMcCann}.

\medskip

We are now ready to introduce the conditions:

\begin{enumerate}
\item[{\bf As}] The cost-sectional curvature is uniformly positive i.e. there exists $C_0 >0$ such that for all $(x,y)\in (\Omega\times \Omega')$, for all $(\nu,\xi) \in \Rd\times\Rd)$ with $\xi\perp\nu$, 
\be
\mathfrak{S_c}(x,y)(\xi,\nu) \geq C_0|\xi|^2|\nu|^2.
\en

\item[{\bf Aw}] The cost-sectional curvature is non-negative: {\bf As} is satisfied with $C_0=0$.

\end{enumerate}

\indent\textsc{Remark on the symmetry of the conditions on $c$.} Let $c^*(y,x):=c(x,y)$, from Proposition \ref{propnewR}, 
one checks that if $c$ satisfies {\bf Aw} (resp. {\bf As}) then $c^*$ satisfies {\bf Aw} (resp. {\bf As} with a different constant). The conditions {\bf A0} and {\bf A2} are also clearly satisfied by $c^*$ if satisfied by $c$.

\subsection{The Riemannian case}
\label{riemannsec}

The construction of optimal maps has been extended in a natural way to smooth compact Riemannian manifolds by McCann in \cite{Mc2} for Lipschitz semi-concave costs. All the above definitions can be translated unambiguously in the Riemannian setting. In particular, the notions of c-exponential map, c-convexity are intrinsic notions (see the Remark 2 after Definition \ref{defiT}).
The definition of cost-sectional curvature \ref{defiR} extends also naturally to the Riemannian setting.  Since it has  been proved in Proposition \ref{propnewR} that the value of the cost-sectional curvature is coordinate-independent, this gives sense to  conditions {\bf Aw}, {\bf As} on a Riemannian manifold.
However, one needs to restrict to the set of pairs $(x,y)$ such that $c$ is smooth in a neighborhood of $(x,y)$, and this becomes an issue for  costs that are functions of the distance: Indeed,  on a compact manifold, the distance can not be smooth on the whole of $M\times M$ (due to the cut-locus). Hence the Riemannian case requires to weaken somehow assumption {\bf A0}.
For  $x$ in $M$, we let $\text{Dom}_x$ be the set of $y$ such that $c(x,y)$ is smooth at $(x,y)$. As developed by the author in \cite{Lsphere}, and with P. Delano\"e in \cite{del-loep}, but also by Y. Kim and R. McCann \cite{KimMcCann}, or by C. Villani in \cite{Vi2}, the relevant geometric condition on $M$ that replaces {\bf A0} is the following:
 {\it for all $ x\in M$, $\mT_x^{-1}(\text{Dom}_x)=-\nabla_xc(x,\text{Dom}_x)$ is convex}.

A case of interest is when  $c(\cdot,\cdot)=\demi d^2(\cdot,\cdot)$ with $d(\cdot,\cdot)$ the distance function (quadratic cost). Then,  the c-exponential map is the exponential map, the map $G_\phi$ will be $x\to\exp_x(\nabla_g\phi)$, the gradient $\nabla_g\phi$ being relative to the Riemannian metric $g$.
(We remind that gradient mappings were first introduced  by X. Cabr\'e \cite{Cabre}, to generalize the Alexandrov-Bakelman-Pucci estimate on Riemannian manifolds.) Then, for  $x$ in $M$, we have $\text{Dom}_x=M\setminus \text{cut-locus}(x)$.
In \cite{Lsphere}, we address the problem of the quadratic cost on the sphere,
as well as the cost $c(x,y)=-\log(|x-y|)$, that appears in the design optimal reflector antenna.
To establish our regularity results, we need to show  \textit{a-priori} that $T(x)$ remains uniformly far from the boundary of $\text{Dom}_x$. This is precisely the object of \cite{del-loep}.

\subsection{Previous regularity results for optimal maps}
The  regularity of optimal maps follows from the regularity of the c-convex potential solution of the Monge-Amp\`ere equation (\ref{ma0}), the former being as smooth as the gradient of the latter. It falls thus  into the theory of viscosity solutions of fully non-linear elliptic equations  \cite{Ca7}, however, the Monge-Amp\`ere equation is degenerate elliptic. 
A very complete reference concerning the regularity theory for the quadratic case are the lecture notes by John Urbas \cite{Urbas}.
Two types of regularity results are usually sought for this type of equations:

{\bf Classical regularity:} show that the equation has classical $C^2$ solutions, provided the measures are smooth enough, and assuming some boundary conditions. Due to the log-concavity of the Monge-Amp\`ere operator, and using  classical elliptic theory (see for instance \cite{GT}),   $C^\infty$ regularity of the solution of (\ref{ma0}) follows from $C^2$ a priori estimates. 

{\bf Partial regularity:} show that a weak solution of (\ref{ma0}) is $C^1$ or $C^{1,\alpha}$ under suitable conditions. We mention also that $W^{2,p}$ regularity results can be obtained.

\paragraph{The Euclidean Monge-Amp\`ere equation and the quadratic cost} \label{background}
 
This corresponds to the case where the cost function is the Euclidean distance squared $c(x,y)=|x-y|^2$ (or equivalently $c(x,y)=-x\cdot y$), for which c-convexity means convexity in the usual sense,  $G_\phi(x)=\nabla\phi(x)$, and  equation (\ref{ma0}) takes the following form
\beq\label{maeucli}
\det D^2 \phi = \frac{\rho_0(x)}{\rho_1(\nabla\phi(x))}.
\enq
Here again, we have $\rho_i={\rm d}\mu_0/{\rm dVol}$, $i=0,1$.
Classical regularity has been established by Caffarelli \cite{Ca1, Ca3, Ca4, Ca5}, Delano\"e \cite{De} and Urbas \cite{U}. The optimal classical regularity result, found in \cite{Ca1, Ca5},  is that for $C^\alpha$ smooth positive densities, and uniformly strictly convex domains, the solution of (\ref{maeucli}) is $C^{2,\alpha}(\bar\Omega)$.
Partial regularity results have been obtained by  Caffarelli {\cite{Ca0, Ca2, Ca3, Ca4}, where it is  shown that for $\mu_0, \mu_1$ having densities bounded away from 0 and infinity, the solution of (\ref{maeucli}) is $C^{1,\alpha}$. Thanks to counterexamples by Wang \cite{Wa2} those results are close to optimal.

\paragraph{The reflector antenna}  The design of reflector antennas can be formulated as a problem of optimal transportation on the unit sphere with cost equal to $-\log|x-y|$.  The potential (height function) $\phi:\Sn^{n-1}\to \R^+$ parametrizes the antenna $A$ as follows:
$A=\{x \phi(x), x\in \Sn^{n-1}\}$. Then the antenna is admissible if and only if $\phi$ is c-convex on $\Sn^{n-1}$ for $c(x,y) = -\log|x-y|$, and $G_\phi(x)$ yields the direction in which the ray coming in the direction $x$ is reflected.
This is the first non quadratic cost for which regularity of solutions has been established. Wang \cite{Wa3, Wa4} (see also Guan and Wang \cite{GuWa} where the results are extended to higher dimension)  has shown classical $C^2$ (and hence $C^\infty$) regularity of solutions of the associated Monge-Amp\`ere equation when the densities are smooth. In a recent work, with totally different techniques, Caffarelli, Huang and Gutierrez \cite{CaGuHu}  have shown  $C^1$ regularity for the solution (i.e. continuity of the optimal map) under the condition that the measures $\mu_0$ and $\mu_1$ have densities bounded away from 0 and infinity. This case of application will also be addressed by our forthcoming paper \cite{Lsphere}. 

\paragraph{General costs and the conditions {\bf As}, {\bf Aw}} Recently an important step was achieved in two papers by Ma, Trudinger, and Wang . They gave in the first paper \cite{MTW} a sufficient condition ({\bf As}, called A3 in their paper) for $C^2$ (and subsequently $C^\infty$) interior regularity. In the second paper \cite{TW}, they could lower this condition down to {\bf Aw} (condition A3w in their paper) to obtain a sufficient condition for global $C^2$ (and subsequently $C^\infty$) regularity, assuming uniform strict c-convexity and smoothness of the domains. Note that  the  result under {\bf Aw} recovers the results of Urbas and Delano\"e  for the quadratic cost. We mention that the results obtained in \cite{MTW, TW} have been exposed by Trudinger in \cite{T}.

\begin{theo}[\cite{TW, T}]
Let $\Omega, \Omega'$ be two  bounded domains of $\Rd$. 
Assume that $\Omega, \Omega'$ are strictly uniformly  c,c*-convex with respect to each other
Let $c, c^*$ satisfy {\bf A0}-{\bf A1}-{\bf A2} and {\bf Aw}  on $\Omega\times\Omega'$. Let $\mu_0, \mu_1$ be two probability measures on $\Omega, \Omega'$ having  densities $\rho_0, \rho_1$. Assume that $\rho_0 \in C^2(\bar\Omega)$ is bounded away from 0, $\rho_1 \in C^2(\bar\Omega')$ is bounded away from 0. 
Then,  for $\phi$ c-convex on $\Omega$ such that $G_{\phi\,\#}\mu_0=\mu_1$, $\phi \in C^3(\Omega)\cap C^2(\bar\Omega)$.
\end{theo}

\noindent
We also mention the continuity result obtained in \cite{McGaShape} concerning optimal transportation between boundaries of uniformly convex domains, that might have some connections with the present work.

\section{Results}

We present some answers to the following four questions:
\begin{enumerate}
\item Is there a sharp necessary and sufficient condition on the cost function  which would guarantee that when both measures have $C^\infty$ smooth densities, and their supports satisfy usual convexity assumptions, the solution of (\ref{ma0}) ( and hence the  optimal map) is $C^\infty$ smooth ? 

\item Is there a necessary and sufficient condition on the cost function and on the data under which optimal maps are continuous ?

\item What are the cost-functions for which  connectedness of the contact set holds (\ref{C}) ?

\item When the cost is set to be the squared distance of a  Riemannian manifold,  what is the meaning of conditions {\bf Aw}, {\bf As} in terms of the Riemannian metric ?

\end{enumerate}

\subsection{Condition {\bf Aw}, connectedness of the contact set and regularity issues }
{\it Answer to questions 1 and 3: Condition  {\bf Aw} is necessary and sufficient  for regularity of optimal maps. Moreover  {\bf Aw} is equivalent to the  connectedness of the contact set.}
\bigskip

In the following theorem, ``smooth'' means $C^\infty$-smooth. This is for simplicity, and one can lower the smoothness assumptions on the domains and the measures, see \cite{TW}.

\begin{theo}\label{main2}
Let $\Omega, \Omega'$ be two  bounded domains of $\Rd$.  Let $c$ be a cost function that satisfies {\bf A0}, {\bf A1}, {\bf A2} on $(\Omega\times\Omega')$. Assume that $\Omega, \Omega'$ are smooth, uniformly strictly c-convex (resp. c*-convex)  with respect to  each other.
The following assertions are equivalent.
\begin{enumerate}
\item The cost function $c$ satisfies {\bf Aw} in $\Omega\times\Omega'$.

\item For  $\mu_0,\mu_1$ smooth strictly positive probability measures in $\bar\Omega, \bar\Omega'$ there exists a  c-convex potential $\phi\in C^1(\Omega)$ such that $G_{\phi\,\#}\mu_0=\mu_1$.
\item For  $\mu_0,\mu_1$ smooth strictly positive probability measures in $\bar\Omega, \bar\Omega'$ there exists a  c-convex potential $\phi\in C^\infty(\bar\Omega)$ such that $G_{\phi\,\#}\mu_0=\mu_1$.

\item For all $\phi$ c-convex in $\Omega$,  for all $x\in \Omega$,  $\partial^c\phi(x)=\partial\phi(x)$.

\item For all $\phi$ c-convex in $\Omega$, for all $x\in \Omega$,  the set  $\{y:\phi(x)+\phi^{c}(y)=-c(x,y)\}$ is c-convex with respect to $x$.

\item Continuously differentiable c-convex potentials are dense among c-convex potentials for the topology of local uniform convergence.
\end{enumerate}
Hence, if condition {\bf Aw} is violated at some points $(x_0,y_0)\in (\Omega\times\Omega')$, there exist smooth positive measures $\mu_0, \mu_1$ on $\Omega, \Omega'$ such that there exists no $C^1$ c-convex potential satisfying $G_{\phi\,\#}\mu_0=\mu_1$.
\end{theo}

\textsc{Remark.} Setting $c^*(y,x)=c(x,y)$ we have seen that $\mathfrak{S_c} \geq 0$ implies $\mathfrak{S_{c^*}} \geq 0$.
Hence all of those assertions are equivalent to their dual counterpart.

We can add the following equivalent condition for {\bf Aw}:
\begin{theo}\label{anoune} Under the assumptions of Theorem \ref{main2}, condition {\bf Aw} holds if and only if, for any $x_0\in \Omega$, $(y_0, y_1) \in \Omega'$, letting $\bar\phi$ be defined by $$\bar\phi(x) = \max\{-c(x,y_0)+c(x_0,y_0), -c(x,y_1)+c(x_0,y_1)\},$$
for any $y_\theta\in [y_0,y_1]_{x_0}$ (see Definition \ref{def-c-seg}),
$$\bar\phi(x) \geq -c(x,y_\theta)+c(x_0,y_\theta)$$
holds in $\Omega$.
\end{theo}
In other words, $f_\theta(x)=-c(x,y_\theta)+c(x_0,y_\theta)$ which is the supporting function that interpolates  at $x_0$ (nonlinearly) between $f_0(x)=-c(x,y_0)+c(x_0,y_0)$ and $f_1(x)=-c(x,y_1)+c(x_0,y_1)$, has to remain below $\max\{f_0,f_1\}$.

\textsc{Remark 1.} The function $\bar\phi$ furnishes the counter-example to regularity when {\bf Aw} is not satisfied, since for a suitable choice of $x_0, y_0, y_1$ $\bar\phi$ can not be approximated by $C^1$ c-convex potentials.

\textsc{Remark 2.} As shown by Propositions \ref{long}, \ref{reciplong}, a quantitative version of Theorem \ref{anoune} holds to express condition {\bf As}.

\textsc{Remark 3.} The assertions Point 1 $\implies$ Points 2, 3, 6 belong to Trudinger and Wang in \cite{TW}.  We show here that condition {\bf Aw} is  necessary: if it is violated at some point, one can always build a counterexample where the solution to (\ref{ma0}) is not $C^1$ even with $C^\infty$ smooth positive measures and good boundary conditions (hence the optimal map is not continuous). Moreover condition {\bf Aw} is equivalent to  a very natural geometric property of c-convex functions.

\subsection{Improved partial regularity under {\bf As}} 

{\it Partial answer to question 2:  There is partial (i.e. $C^1$ and $C^{1,\alpha}$) regularity under {\bf As}, requiring much lower assumptions on the measures than what is needed in the quadratic case. There can not be $C^1$ regularity without {\bf Aw}. When only {\bf Aw} is satisfied, the question of $C^1$ regularity remains open, except for the case $c(x,y)=|x-y|^2$ treated by Caffarelli \cite{Ca3}.}

\bigskip

Let us begin by giving the two integrability conditions that will be used in this result.
The first one reads
\beq\label{Lp}
&&\nonumber\text{ For some } p\in ]n, +\infty], C_{\mu_0}>0,\\
&& \mu_0(B_\epsilon(x)) \leq C_{\mu_0} \epsilon^{n (1-\frac{1}{p})} \text{ for all } \epsilon \geq 0, x \in \Omega.
\enq  
The second condition reads
\beq\label{mini}
&& \nonumber\text{ For some } f:\R^+\to \R^+ \text{ with }\lim_{\epsilon \to 0} f(\epsilon) = 0, \\
&&\mu_0(B_\epsilon(x)) \leq f(\epsilon) \epsilon^{n (1-\frac{1}{n})} \text{ for all } \epsilon \geq 0, x \in \Omega.
\enq
In order to appreciate the forthcoming theorem, let us mention a few facts on these integrability conditions (the proof of this proposition is given at the end of the paper). 
\begin{prop}\label{proprho}Let $\mu_0$ be a probability measure on $\Rd$.
\begin{enumerate}
\item If $\mu_0$ satisfies (\ref{Lp}) for some $p>n$, $\mu_0$ satisfies (\ref{mini}).
\item If $\mu_0 \in L^p(\Omega)$ for some $p>n$, $\mu_0$ satisfies (\ref{Lp}) with the same $p$.
\item If $\mu_0 \in L^n(\Omega)$, $\mu_0$ satisfies (\ref{mini}).
\item If $\mu_0$ satisfies (\ref{mini}), $\mu_0$ does not give mass to set of Hausdorff dimension less than or equal to $n-1$, hence (\ref{mini}) guarantees the existence of an optimal map.
\item There are probability measures on $\Omega$ that satisfy (\ref{Lp}) (and hence (\ref{mini})) and that are not absolutely continuous with respect to the Lebesgue measure.
\end{enumerate}
\end{prop}
Then our result is

\begin{theo}\label{main}
Let $c$ be a cost function that satisfies assumptions {\bf A0}, {\bf A1}, {\bf A2}, {\bf As} on $(\Omega\times\Omega')$, $\Omega, \Omega'$ being bounded domains uniformly strictly c(resp c*)-convex with respect to each other. 
Let $\mu_0, \mu_1$ be probability measures  respectively on $\Omega$ and $\omega'\subset\Omega'$, with $\omega'$ c-convex with respect to $\Omega$.   Let $\phi$ be a c-convex potential on $\Omega$ such that $G_{\phi\,\#}\mu_0=\mu_1$.
Assume that $\mu_1  \geq m\,{\rm dVol}$ on $\omega'$ for some $m>0$.
\begin{enumerate}
\item  Assume that $\mu_0$ satisfies (\ref{Lp}) for some $p>n$.
Let $\alpha = 1-\frac{n}{p},  \  \beta = \frac{\alpha}{4n-2+\alpha}$.
Then for any $\delta>0$ we have
\be
\|\phi \|_{C^{1,\beta}(\Omega_\delta)} \leq {\cal C},
\en 
and ${\cal C}$ depends only 
on $\delta>0$, $C_{\mu_0}$ in (\ref{Lp}), on $m$,  on the constants  in conditions {\bf A0}, {\bf A1}, {\bf A2}, {\bf As} and on $C_\mT$ in (\ref{bma}).

\item If $\mu_0$ satisfies (\ref{mini}),
then $\phi$ belongs to $C^1(\Omega_\delta)$ and the modulus of continuity of $\nabla\phi$ is controlled  by $f$ in (\ref{mini}).
\end{enumerate}
\end{theo}
As an easy corollary of Theorem \ref{main}, we can extend the $C^1$ estimates to the boundary
if the support of the measure $\mu_0$ is compactly contained in $\Omega$.
\begin{theo}\label{main4}
Assume in addition to the assumptions of Theorem \ref{main} that $\mu_0$ is supported in $\bar\omega$, with $\omega$ compactly contained in $\Omega$. Then, if $\mu_0$ satisfies (\ref{mini}), $\phi \in C^1(\bar\omega)$ and if $\mu_0 $ satisfies (\ref{Lp}), $\phi \in C^{1,\beta}(\bar\omega)$, with $\beta$ as in Theorem \ref{main}.
\end{theo}

\textsc{Remark on the conditions on $\Omega, \Omega'$.} 
Our result holds true for $\mu_0$ supported in any subset $\omega$ of $\Omega$ (hence not necessarily c*-convex), and $\mu_1$ supported in any subset $\omega'$ of $\Omega'$ c-convex (but not necessarily strictly) with respect to $\Omega$.
Hence what we need is the existence of supersets $\Omega, \Omega'$ uniformly c(c*)-convex with respect to each other, in order  to use the results of \cite{TW}. The only point where we need this condition is during the proof of Proposition \ref{supporting}, where we rely on Theorem \ref{main2} to assert $\partial\phi = \partial^c\phi$.
However, in \cite{MTW}, Ma, Trudinger and Wang proved the following:
\begin{theo}[\cite{MTW}]\label{theoTW2}
Let $c$ satisfy {\bf A0}, {\bf A1}, {\bf A2},  {\bf As}, on $\Omega \times\Omega'$,
 $\Omega'$ being c-convex with respect to $\Omega$.
Then, for $\mu_0$, $\mu_1$ $C^2$ smooth positive probability measures on $\Omega, \Omega'$, the c-convex potential $\phi$ such that $G_{\phi\,\#}\mu_0=\mu_1$ is $C^2$ smooth inside $\Omega$.
\end{theo}
Using this result,  Proposition \ref{egal} yields that for all $\phi$ c-convex on $\Omega$, $\partial^c\phi = \partial\phi$.
Hence we could have relaxed the assumptions of Theorem \ref{main} on $\Omega, \Omega'$,  only requiring $\Omega'$ to be c-convex with respect to $\Omega$, (i.e. no c*-convexity on $\Omega$, no strict c-convexity of $\Omega'$).
Note that the proof of Theorem \ref{theoTW2} has been completed later on by Trudinger and Wang in \cite{TW2}, relying in part on our Proposition \ref{long} (which is an independent result).
Thus we can now state the following result:
\begin{theo}
The results of Theorem \ref{main} hold assuming only for $\Omega, \Omega'$ that $\Omega'$ is c-convex  
with respect to $\Omega$.
\end{theo}
We mention  that  the results of Kim and McCann \cite{KimMcCann}, obtained simultaneously with those of \cite{TW2} but using different techniques, allow also to  complete the proof of
Theorem \ref{theoTW2}, under the assumption that $\Omega'$ and $\Omega$ are c-(c*)-convex with respect to each other.
This allows to drop the strict convexity assumption in Theorem \ref{main}. 

\textsc{Remark on the integrability conditions.}  The integrability conditions on $\mu_0, \mu_1$ are really mild: we only ask that $\mu_1$ be bounded by below, and that $\mu_0(B_r) \leq r^{n-p}$ for $p\geq 1$ ($p >1$ yields $C^{1,\alpha}$ regularity) (see conditions (\ref{Lp}) and (\ref{mini})  and the subsequent discussion). The continuity of the optimal map is also asserted in the case $\mu_0 \in L^n$ (that implies (\ref{mini})), which is somehow surprising: indeed $D^2\phi \in L^n$ does not imply $\phi \in C^1$, but here $\det(D^2\phi-{\cal A}(x,\nabla\phi))\in L^n$ implies $\phi \in C^1$.
In a forthcoming work, we shall show that our result adapts to the reflector antenna, hence improving the result  obtained independently by Caffarelli,  Gutierrez and Huang \cite{CaGuHu} on reflector antennas.   Moreover our techniques yield quantitative $C^{1,\alpha}$ estimates: the exponent $\alpha$ can explicitly computed. Finally, our continuity estimates extends up to the boundary (Theorem \ref{main4}). This is achieved through a geometric formulation of condition {\bf As}.

A full satisfactory answer would include a general result of partial regularity under condition {\bf Aw}. This result is expected in view of the Euclidean case (since the quadratic cost is really the limit case for condition {\bf Aw}). Note that, in view of counterexamples  given in \cite{Wa2}, the results under {\bf Aw} can not be as good as under {\bf As}, and can not be much better than Caffarelli's results \cite{Ca3} that require densities bounded away from 0 and infinity.

\subsection{Conditions {\bf Aw}, {\bf As} for the quadratic cost of a Riemannian manifold }
We refer the reader to the remark 2. after the definition of the cost-sectional curvature (\ref{defR}) where the intrinsic meaning of (\ref{defR}) on a manifold is discussed.

{\it Partial answer to question 4:  When the cost is the Riemannian distance squared, the cost-sectional curvature at $y=x$ equals (up a multiplicative constant) the Riemannian sectional curvature}

\begin{theo}\label{costcurvature} Let $M$ be a $C^4$ Riemannian manifold. Let $c(x,y)=d^2(x,y)/2$ for all $(x,y)\in M\times M$. Let $\mathfrak{S_c}$ be given by (\ref{defR}), 
Then, for all $\xi,\nu \in T_xM$, 
\be
\frac{\mathfrak{S_c}(x,x)(\nu,\xi)}{|\xi|^2_g|\nu|^2_g-(\xi\cdot\nu)_g^2}=\frac{2}{3}\cdot\text{ Sectional Curvature of } M \text{ at }x \text{ in the 2-plane }(\xi,\nu).
\en
Hence if {\bf Aw} (resp, {\bf As}) is satisfied at $(x,x)$,  the sectional curvature of $M$ at $x$  is non-negative (resp. strictly positive).
\end{theo}

\begin{cor}
Let $M$ be a compact Riemannian manifold. If the sectional curvature of $M$ is not everywhere non-negative, there are smooth positive measures on $M$ such that the optimal map (for the cost function $c(x,y) =d^2(x,y)/2$) is not continuous.
\end{cor} 
At the end of the proof of Theorem \ref{costcurvature}, we give a counterexample to regularity for a two-dimensional manifold with negative sectional curvature.

This observation closes (with a negative answer) the open problem of the regularity of optimal gradient maps when the manifold does not have non-negative sectional curvature everywhere.
There is a partial converse assertion in the special case of constant sectional curvature:
{\it The quadratic cost on the round sphere $\Sn^{n-1}$ satisfies {\bf As}}. This will be the object of a forthcoming work \cite{Lsphere}. Hence our previous result can be adapted to this Riemannian case.

\subsection{Examples of costs that satisfy {\bf As} or {\bf Aw}}

We repeat the collection of cost that was given in \cite{MTW}, and \cite{TW}.

\begin{itemize}
\item $c(x,y)= \sqrt{1+|x-y|^2}$ satisfies {\bf As}.
\item $c(x,y)= \sqrt{1-|x-y|^2}$ satisfies {\bf As}.
\item $c(x,y)=(1+|x-y|^2)^{p/2}$ satisfies {\bf As} for $1\leq p<2$, $|x-y|^2 < \frac{1}{p-1}$. 
\item $c(x,y)=|x-y|^2+ |f(x)-g(y)|^2$ $f,g:C^4(\Rd;\R)$   convex  (resp.  strictly convex) with $|\nabla f|, |\nabla g| < 1$ satisfies  {\bf Aw} (resp. {\bf As}).
\item $c(x,y) = \pm\frac{1}{p}|x-y|^p, p\neq 0$ and  satisfies {\bf Aw} for $p=\pm 2$ and {\bf As} for $-2 < p < 1$ ($-$ only). 
\item $c(x,y) =-\log|x-y|$ satisfies {\bf As} on $\Rd\times\Rd \setminus \{(x,x), x\in \Rd\}$.
\item The reflector antenna problem  (\cite{Wa3}) corresponds to the case $c(x,y)=-\log|x-y|$ restricted to $\Sn^n$. As pointed out in \cite{TW}, this cost satisfies {\bf As} on $\Sn^{n-1}\times \Sn^{n-1}\setminus \{x=y\}$. 
\item As shown in a forthcoming paper \cite{Lsphere}, the squared Riemannian distance on the sphere satisfies {\bf As} on the set $\Sn^{n-1}\times \Sn^{n-1}\setminus \{x=-y\}$. Note that it is the restriction to $\Sn^{n-1}$ of the cost $c(x,y)=\theta^2(x,y)$, where $\theta$ is the angle formed by $x$ and $y$. 
(For those two cases, see paragraph \ref{riemannsec} where the meaning of conditions {\bf Aw}, {\bf As} on a Riemannian manifold is discussed).

\end{itemize}

\paragraph{Acknowledgments}
At this point I wish to express my gratitude to Neil Trudinger and Xu-Jia Wang, for many discussions, and for sharing results in progress while we were all working on this subject. I also thank C\'edric Villani for fruitful discussions, and Philippe Delano\"e with whom we started to think about the problem of regularity for optimal transportation on the sphere. I also wish to thank Robert McCann,
who first raised to me the issue of the connectedness of the contact set, in 2003. 
I gratefully acknowledge the support of a French Australian exchange grant PHC FAST EGIDE No.12739WA. I am also grateful for the hospitality of the Center for Mathematics and its Applications at University of Canberra.

\bigskip

\section{Proof of Theorem \ref{main2}}

We begin with the following uniqueness result of independent interest:
\begin{prop}\label{unicite}
Let $\mu, \nu$ be two probability measures on $\Omega, \Omega'$, with $\Omega$ and $\Omega'$ connected domains of $\Rd$. Assume that either $\mu$ or $\nu$ is positive Lebesgue almost everywhere in $\Omega$ (resp. in $\Omega'$). Then, among all pairs of functions $(\phi, \psi)$ such that $\phi$ is c-convex, $\psi$ is c*-convex, the problem (\ref{kanto}) has at most one minimizer up to an additive constant.
\end{prop}

The proof of this proposition is deferred to the end of the paper.

\subsection{Condition {\bf Aw} implies connectedness of the contact set}

We will begin with the following lemma:
\begin{lemme}\label{grad-converge}
Let $\phi$ be  c-convex. Let $(\phi\epu)_{\epsilon>0}$ be a sequence of  c-convex potentials that converges uniformly to $\phi$ on compact sets of $\Omega$.  Then, if $p=-\nabla_xc(x_0,y)\in \partial\phi(x_0)$, $x_0\in \Omega, y\in \Omega'$, there exists a sequence $(x\epu)_{\epsilon>0}$ that converges to $x_0$, a sequence $(y\epu)_{\epsilon>0}$ that converges to $y$ such that  $p\epu=-\nabla_xc(x\epu,y\epu)\in\partial\phi\epu(x\epu)$.
Finally, $p\epu$  converges to $p$. 
\end{lemme}

\textsc{Proof.} Let $y=\mathfrak{T}_{x_0}(p)$, i.e. $p=-\nabla_xc(x_0,y)$.
Since $\phi, \phi\epu$ are c-convex and $c$ is semi-concave, there exists $K,r>0$  so that 
\be
&&\tilde\phi(x):=\phi(x)+K|x-x_0|^2/2+c(x,y),\\
&& \tilde\phi\epu(x):=\phi\epu(x)+K|x-x_0|^2/2+c(x,y),
\en
are convex on $B_r(x_0)$ compactly contained in $\Omega$.
One can also assume, by subtracting a constant  that $\tilde\phi(x_0)=0$, and that $\tilde\phi(x)\geq 0$ on $\Omega$. 
Finally, one can assume (by relabeling the sequence) that on $B_r(x_0)$  we have $|\phi\epu - \phi| \leq \epsilon$.

Consider then $\tilde\phi\epu^\delta=\tilde\phi\epu +\delta |x-x_0|^2/2-\epsilon$.
We have $\tilde\phi\epu^\delta(x_0) \leq 0$, and on $\partial B_\mu(x_0)$, with $\mu\leq r$,  
\be
\tilde\phi\epu^\delta(z)&\geq& \tilde\phi(z)+\delta \mu^2/2 - 2\epsilon\\
&\geq& \delta \mu^2/2 -2\epsilon.
\en 
By taking $\mu=\epsilon^{1/3}, \delta=4\epsilon^{1/3}$, we get that $\tilde\phi\epu^\delta$ has a local minimum in $B_\mu(x_0)$, hence at some point $x\epu\in B_\mu(x_0)$, we have 
\be
\partial\phi\epu(x\epu)\owns -\nabla_x c(x\epu,y)-K(x\epu-x_0) - \delta (x\epu-x_0).
\en
Then we have $|(K+\delta)(x\epu-x_0)|$ small, and thanks to {\bf A1}, {\bf A2}, there exists $y\epu$ close to $y$ such that $\nabla_xc(x\epu,y\epu) = \nabla_x c(x\epu,y) + K(x\epu-x_0)+ \delta (x\epu-x_0)$. Thus $-\phi\epu(x) - c(x,y\epu)$ has a critical point at $x\epu$. This implies that $p\epu = -\nabla_xc(x\epu,y\epu) \in \partial\phi\epu(x\epu)$. Finally, since $x\epu\to x, y\epu\to y$, we conclude $p\epu \to p$.  

$\hfill\Box$
\\

Now we prove that $\partial^c\phi=\partial\phi$. In order to do this, we must show that if $\phi$ is c-convex, if $-\phi(\cdot)-c(\cdot, y)$ has a critical point at $x_0$, this is a global maximum.

We first have the following observation:
\begin{lemme}\label{nonC1}
Let $\phi$ be c-convex. Assume that $-\phi - c(\cdot,y)$ has a critical point at $x_0$ (i.e. $0\in \partial\phi(x_0)+\nabla_xc(x_0,y)$), and that it is not a global maximum.  Then $\phi$ is not differentiable at $x_0$.
\end{lemme}

\textsc{Proof.} Indeed,  $-\phi(\cdot) - c(\cdot,y)$ has a critical point at $x_0$, but we don't have $\phi(x_0) + \phi^c(y)=-c(x_0,y)$. However, there is a point $y'$ such that $\phi(x_0) + \phi^c(y')=-c(x_0,y')$. Hence, $\{ -\nabla_x c(x_0,y), -\nabla_x c(x_0,y')\} \in \partial\phi(x_0)$, and we have $\nabla_xc(x_0,y)\neq\nabla_xc(x_0,y')$ from assumption {\bf A1}. 

$\hfill \Box$

We show the following:
\begin{prop}\label{egal}
Assume {\bf D} holds.
Let $p=-\nabla_xc(x_0,y)\in \partial\phi(x_0)$ with $\phi$ c-convex.
Then  $-\phi(\cdot) - c(\cdot, y) $ reaches a global maximum at  
$x_0$.
\end{prop}

\begin{itemize}
\item[{\bf D}] $C^1$ c-convex functions are dense in the set $\Big\{\phi$ c-convex on $\Omega$,  $G_\phi(\Omega)\subset\Omega'\Big\}$ for the topology of uniform convergence on compact sets of $\Omega$.
\end{itemize}

\textsc{Proof.} Assume the contrary, i.e. that $-\phi(x_1) - c(x_1,y) > -\phi(x_0)-c(x_0,y)$ for some $x_1\in \Omega$.
We use {\bf D}: there exists a sequence of $C^1$ c-convex potentials $(\phi\epu)_{\epsilon>0}$ that converges to $\phi$. 
 We use Lemma \ref{grad-converge}: there will exist a sequence $(x\epu)_{\epsilon>0}$ such that $x\epu\to x_0$ and $\nabla\phi\epu(x\epu)\to -\nabla_x c(x_0,y)$. Let $y\epu$ be such that $\nabla\phi\epu(x\epu) = -\nabla_xc(x\epu, y\epu)$. Then $y\epu\to y$.
 Since $\phi\epu$ is $C^1$, by Lemma \ref{nonC1}, $x\epu$, the critical point of $-\phi\epu(\cdot) - c(\cdot,y\epu)$ is necessarily a global maximum.  Finally, since $\phi\epu$ converges locally uniformly to $\phi$, we see that $-\phi(\cdot) -c(\cdot,y)$ reaches at $x_0$ a  global maximum.

\begin{lemme}\label{A3D}
Assume $\Omega, \Omega'$ are bounded, uniformly strictly c-(c*-)  convex with respect to  each other.
Assume that $c$ satisfies {\bf A0}, {\bf A1}, {\bf A2}, {\bf Aw} on $\Omega\times \Omega'$.   Then {\bf D} holds.
\end{lemme}

\textsc{Proof.} As we will see, this result is implied immediately by the result of \cite{TW}  combined with  Proposition \ref{unicite}.
Let $\phi$ be c-convex. Denote $\mu_1 = G_{\phi\, \#}{\bf 1}_{\Omega}{\rm dVol}$. 
Note that from Proposition \ref{unicite}, $\phi$ is the unique up to a constant c-convex potential such that $G_{\phi\, \#}{\bf 1}_{\Omega}{\rm dVol} =\mu_1$.
Consider a sequence of smooth positive densities $(\mu_1\ep)_{\epsilon > 0}$ in $\Omega'$ such that $\mu_1\ep{\rm dVol}$  converges weakly-$*$ to $\mu_1$, and has same total mass than $\mu_1$. Consider $\phi\epu$ such that $G_{\phi\epu\,\#}{\bf 1}_{\Omega}{\rm dVol} = \mu_1\ep {\rm dVol}$. From \cite{TW}, $\phi\epu$ is $C^2$ smooth inside $\Omega$. 
Then, by Proposition \ref{unicite}, up to a normalizing constant, $\phi\epu$ is converging to $\phi$, and $\nabla\phi\epu$ is converging to $\nabla\phi$ on the points where $\phi$ is differentiable. 
$\hfill\Box$ 

Hence, under the assumptions of Lemma \ref{A3D}, $\partial\phi(x)=\partial^c\phi(x)$.
In view of Proposition \ref{remarks},  the equality $\partial\phi(x)=\partial^c\phi(x)$ for all $\phi, x$ is equivalent to the c-convexity of the set 
\be
G_\phi(x) = \Big\{y: \phi(x)+\phi^c(y)=-c(x,y)\Big\}.
\en
This shows that condition {\bf Aw} is sufficient. 
$\hfill\Box$

\subsection{Condition {\bf Aw} is necessary for smoothness and connectedness of the contact set}

We now show that if {\bf Aw} is violated somewhere in $(\Omega\times\Omega')$, there will exist a c-convex potential for which we don't have $\partial\phi=\partial^c\phi$. Assuming this, in view of Lemma \ref{A3D} and Proposition \ref{egal}, this will imply that this potential can not be a limit of $C^1$-smooth c-convex potentials. Hence, considering the sequence $(\phi\epu)_{\epsilon > 0}$ used  in the proof of Lemma \ref{A3D}, this sequence will not be $C^1$ for $\epsilon$ smaller than some $\epsilon_0$.  This implies in turn that there exists smooth positive densities $\mu_0, \mu_1$ in $\Omega, \Omega'$ such that the c-convex potential $\phi$ satisfying $G_{\phi\,\#}\mu_0=\mu_1$ is not $C^1$ smooth.

Assume that for some $x_0 \in \Omega, y\in \Omega', p=-\nabla_xc(x_0,y)$, for some $\xi,\nu$ unit vectors in $\Rd$ with $\xi\perp\nu$, one has
\beq\label{pasA3}
D^2_{p_\nu p_\nu}\left[p\to D^2_{x_\xi x_\xi}c(x,\mathfrak{T}_x(p))\right] \geq N_0>0.
\enq
Let $y_0=\mathfrak{T}_{x_0}(p-\epsilon \nu), y_1=T_{x_0}(p+\epsilon \nu)$, with $\epsilon$ small, and recall that $y=\mathfrak{T}_{x_0}(p)$.
Hence $y$ is the 'middle' of the c-segment $[y_0,y_1]_x$.
Let us define 
\beq\label{defbarphi}
\bar \phi(x)= \max\Big\{-c(x,y_0) + c(x_0,y_0), -c(x,y_1) + c(x_0,y_1)\Big\}.
\enq
(This function will be used often in the geometric interpretation of {\bf As}, {\bf Aw}. It is the ``second simplest'' c-convex function, as the supremum of two supporting functions. It plays the role of $(x_1,...,x_n)\to |x_1|$ in the Euclidean case.)

Note first that $\xi\perp\nu$ implies that $\xi\perp (\nabla_xc(x_0,y_1)-\nabla_xc(x_0,y_0))$.
Consider near $x_0$ a smooth curve $\gamma(t)$ such that $\gamma(0)=x_0$, $\dot\gamma(0)=\xi$, and such that for $t\in [-\delta, \delta]$, one has 
\be
f_0(\gamma(t)):=-c(\gamma(t),y_0) + c(x_0,y_0)= -c(\gamma(t),y_1) + c(x_0,y_1)=:f_1(\gamma(t)).
\en
Such a curve exists by the implicit function theorem, and it is $C^2$ smooth.
On $\gamma$, we have
\be
\bar\phi = \demi (f_0 + f_1)
\en since $f_0=f_1$ on $\gamma$.
Then we compare $\demi (f_0 + f_1)$ with $-c(x,y)+c(x_0,y)$.
By (\ref{pasA3}) we have
\be
\demi\left[D^2_{x_\xi x_\xi} c(x_0,y_0) + D^2_{x_\xi x_\xi} c(x_0,y_1)\right] \geq D^2_{x_\xi x_\xi} c(x_0,y)+c(\epsilon, N_0),
\en
where $c(\epsilon, N_0)$ is positive for $\epsilon$ small enough.
Then of course $\nabla_x c(x_0,y) = \demi[\nabla_x c(x_0,y_0)+ \nabla_x c(x_0,y_1)]$.
Hence we have, for $\epsilon$ small enough,
\be
&&[-c(\gamma(t),y)+c(x_0,y)]-\bar\phi(\gamma(t))\\
&=&[-c(\gamma(t),y)+c(x_0,y)]- \demi (f_0 + f_1)(\gamma(t))\\
&=&\Big[ \demi\left[D^2_{x x} c(x_0,y_0) + D^2_{x x} c(x_0,y_1)\right]- D^2_{x x} c(x_0,y)\Big] \cdot (\gamma(t)-x_0)\cdot(\gamma(t)-x_0)/2  + o(t^2)\\
&=& \Big[ \demi\left[D^2_{x_\xi x_\xi} c(x_0,y_0) + D^2_{x_\xi x_\xi} c(x_0,y_1)\right]- D^2_{x_\xi x_\xi} c(x_0,y)\Big] t^2/2 + o(t^2)\\
&\geq& c(\epsilon, N_0)t^2/2 + o(t^2).
\en
This will be strictly positive for $t\in [-\delta,\delta]\setminus\{0\}$ small enough, and of course the difference $-\bar \phi - [c(x,y)-c(x_0,y)]$ vanishes at $x_0$. 
Obviously, the function $\bar \phi$ is c-convex, $-\bar\phi(\cdot) - c(\cdot,y)$ has a critical point at $x_0$, and this is not a global maximum. Hence, from Proposition \ref{egal},  {\bf D} can not hold true.

The proof of Theorems \ref{main2}, \ref{anoune} is complete.

$\hfill\Box$.

\section{Proof of Theorem \ref{main}}

\subsection{Sketch of the proof} The key argument of the proof is the geometrical translation of condition {\bf As}: assume that $\phi$   $c$-convex  is not differentiable at $x=0$, hence, for some pair $y_0, y_1$, 
$-\phi(\cdot) - c(\cdot,y_0)$ and $-\phi(\cdot) - c(\cdot,y_1)$  both reach a  maximum at $x=0$. 
(From Theorem \ref{main2}, under  {\bf As}, all critical points of  $-\phi(\cdot) - c(\cdot,y)$ are global maxima.)
Consider  $y_\theta$ in the c-segment with respect to $x=0$ joining $y_0$ to $y_1$. As we will see in Proposition \ref{long},  
the function $-\phi(\cdot) - c(\cdot,y_\theta)$  will  reach a  maximum at $x=0$,
and  condition {\bf As} implies moreover that for $\theta\in [\epsilon, 1-\epsilon]$ this maximum will be strict in the following sense:  we will have 
$$-\phi(x)  -c(x,y_\theta) \leq -\phi(0) - c(0,y_\theta) - \delta |x|^2 + o(|x|^2),$$
with $\delta>0$ depending on $|y_1-y_0|$ and $C_0>0$ in condition {\bf As}, and bounded by below for $\theta$ away from $0$ and $1$.

Then, by estimating all supporting functions to $\phi$ on $B_\eta(0)$ a small ball centered at 0 , we will find that for $y$ in
a $C\eta$ neighborhood of $\{y_\theta\}_{\theta\in [1/4, 3/4]}$, $C>0$ depending on $\delta$ above, $-\phi(\cdot) -c(\cdot, y)$ will reach a local maximum in $B_\eta(0)$. 
 Hence $G_\phi(B_\eta(0))$ contains a $C\eta$ neighborhood of $\{y_\theta\}_{\theta\in [1/4, 3/4]}$. 
This is the Proposition \ref{supporting}. Once this is shown, we can  contradict the bound on the Jacobian determinant of $G_\phi$.

We now enter into the rigorous proof of Theorem \ref{main}, this proof is articulated in three parts.

\subsection{Part I. Geometric interpretation of condition {\bf As}}
 
This proposition is the geometrical translation of assumption {\bf As}. Actually, 
as we will see in Proposition \ref{reciplong}, the result of Proposition \ref{long} is equivalent to assumption {\bf As} for a smooth cost function.

\begin{prop}\label{long}Let $c$ be a cost function that satisfies $\bf{A0}, \bf{A1}, \bf{A2}, \bf{As}$ on $\Omega \times \Omega'$. 
For $x_0\in \Omega$, $y_0, y_1 \in \Omega'$,  let $\{y_\theta\}_{\theta\in [0,1]}$ be the c-segment with respect to $x_0$ joining $y_0$ to $y_1$, in the sense of Definition \ref{def-c-seg}, and assume that $\Omega'$ is c-convex with respect to $x_0$.
Let $$\bar \phi(x)=\max\{-c(x,y_0)+ c(x_0,y_0), -c(x,y_1)+ c(x_0,y_1)\}.$$ There exist constants      $\delta_0, C>0$ and $\gamma$ such that for all $\epsilon\in ]0,\demi[$,  $\theta\in [\epsilon,1-\epsilon]$, for all $x\in\Omega$ such that $|x-x_0| \leq C\epsilon$, 
\be
\bar \phi (x) \geq -c(x,y_\theta) + c(x_0, y_\theta) + \delta_0 \theta(1-\theta) |y_1-y_0|^2 |x-x_0|^2 -\gamma|x-x_0|^3,
\en
with  lower bounds on $\delta_0$ and $C$  and an upper bound on $\gamma$ 
that depend on the bounds in assumptions $\bf{A0}, \bf{A2}, \bf{As}$, on an upper bound on $|y_1-y_0|$,  and on  $C_\mT$ in (\ref{bma}).

\end{prop}

\noindent
\paragraph{Preliminary Results}  Shifting and rotating the coordinates, we can assume that $x_0=0$
and that $\nabla_x c(0, y_0) - \nabla_x c(0, y_1)$ is parallel to $e_1$. 
Then, we observe  the following fact:
\begin{prop}\label{subtract}
Subtracting from $c$ a smooth function $x\to \lambda(x)$ that depends only on $x$ does not change the  map solution of the optimal transportation problem, and the new cost $c(x,y) - \lambda(x)$ will still satisfy assumptions {\bf A0}, {\bf A1}, {\bf A2},  {\bf Aw}. The optimal potential will be changed according to the rule $\phi \to \phi + \lambda$.
If moreover the function $\lambda$ is affine, this modification does not change the bounds in assumptions {\bf A2}, {\bf As}.
\end{prop}

\noindent
Using Proposition \ref{subtract}, we can subtract from $c$ the affine function given by $$\lambda(x) = \nabla_x c(0, y_0)\cdot(x-x^1e_1),$$ so that the new cost $c$ will satisfy
\beq
\nabla_x c(0, y_0) = -a e_1, \;\;\nabla_x c(0, y_1) = -b e_1,\label{defa}
\enq
 for some $a \neq b$ from assumption {\bf A1} (we will assume hereafter that $b>a$).
Note that (\ref{defa}) is equivalent to $$y_0 = \mathfrak{T}_{x=0}(a e_1),\;\; y_1=\mathfrak{T}_{x=0}(b e_1).$$
We then have for all $\theta \in [0,1]$,
\beq
-c(x,y_\theta)+c(0,y_\theta)
&=&[\theta b + (1-\theta)a] x^1 -\demi D^2_{xx}c(0,y_\theta)\cdot (x, x) +o(|x|^2)\label{ft}.
\enq

We now have the following  Lemma, which is the point where we use assumption {\bf As}.
\begin{lemme}\label{deltaDelta}
Under assumptions and notations of Proposition \ref{long}, in particular assuming {\bf As}, for all $x\in \R^n$, for all $\theta \in [0,1]$, letting 
$a$ and $b$ be defined through (\ref{defa}), one has
\be
-D^2_{xx}c(0,y_\theta)\cdot (x,x) &\leq& - \left[ (1-\theta) D^2_{xx}c(0,y_0) + \theta D^2_{xx}c(0,y_1)\right]\cdot (x,x)\\
&& - \delta |x|^2 \\
&& + \Delta|x_1|^2,
\en
where 
\be
\delta &=& \frac{1}{4} C_0|b-a|^2\theta(1-\theta),\\
\Delta &=& \frac{\Delta_0^2}{C_0} |b-a|^2\theta(1-\theta),
\en
with  $C_0$   given in assumption {\bf As} 
and $\Delta_0$ depending on  $\|c(\cdot,\cdot)\|_{C^4(\Omega\times\Omega')}, \|[D_{xy} c]^{-1}\|_{\Linf(\Omega\times\Omega')}.$
Note in particular that under {\bf A0}, {\bf As}, $C_0$ is bounded away from 0 and $+\infty$.

\end{lemme}
We will also need the following elementary estimates, that we state without proof:

\begin{lemme}\label{lemmeC1}
Under assumptions and notations of Proposition \ref{long}, for all $x\in \R^n$, for all $\theta, \theta' \in [0,1]$, 
\beq\label{defC1}
&&\demi|D^2_{xx}c(0,y_\theta)\cdot (x, x) - D^2_{xx}c(0,y_{\theta'})\cdot (x,x)| \leq C_1 |\theta-\theta'||x|^2,
\enq
with $C_1$  depends on $|b-a|$, $\|[D_{xy}c]^{-1}\|_{\Linf(\Omega\times\Omega')}$ and $\|c(\cdot, \cdot)\|_{C^3(\bar\Omega \times \bar\Omega')}$.
\end{lemme}

\begin{lemme}\label{semi-convex}
Let $[t_0,t_1]\subset \R$ and $f$ belong to  $C^2([t_0, t_1], \R)$.
\begin{enumerate}
\item If $f'' \geq \alpha $,  we have, for all $t_0,t_1 \in \R$, 
\be
\theta f(t_0) + (1-\theta)f(t_1) \geq f(\theta t_0 + (1-\theta)t_1) + \demi\alpha\theta (1-\theta) |t_1-t_0|^2. 
\en
\item In all cases we have 
\be
\big|\theta f(t_0) + (1-\theta)f(t_1)- f(\theta t_0 + (1-\theta)t_1)\big| \leq \demi \|f\|_{C^2(t_0, t_1)}\theta (1-\theta) |t_1-t_0|^2. 
\en
\end{enumerate}
\end{lemme}

\bigskip

\textsc{Proof of Lemma \ref{deltaDelta}.}   We apply the first part of Lemma \ref{semi-convex} to the function
\be
f: t \to -D^2_{xx} c(0,\mathfrak{T}_{x=0}(t e_1))\cdot (x', x')
\en
where $x'$ is equal to $(0,x^2,..,x^n)$, and hence $x' \perp e_1$.
From assumption {\bf As}, this function  satisfies $f'' \geq  C_0 |x'|^2$. Then, by choosing $t_0=a, t_1= b$ (note that $y_\theta = \mathfrak{T}_{x=0}((\theta b + (1-\theta) a)e_1)$), we obtain that
\be
-D^2_{xx}c(0,y_\theta)\cdot (x', x') &\leq& - \left[(1-\theta) D^2_{xx}c(0,y_0) + \theta D^2_{xx}c(0,y_1)\right]\cdot (x', x') \\ 
&&- \demi C_0 |x'|^2 \theta(1-\theta)|b-a|^2 .
\en 
To conclude the lemma, we  have to control of the terms where $x^1$ appears.
For this we apply the second part of Lemma \ref{semi-convex} to
\be
g: t \to D^2_{xx} c(x,\mathfrak{T}_x(t e_1))\cdot (x,x) - D^2_{xx} c(x,\mathfrak{T}_x(t e_1))\cdot (x',x'),
\en 
for which we have
$|g''| \leq 2\Delta_1 |x^1||x|$, where $\Delta_1$ depends on $\|c(\cdot,\cdot)\|_{C^4}$ and on $\|[D_{xy}c]^{-1}\|_{\Linf}$. This yields
\be
-D^2_{xx}c(0,y_\theta)\cdot (x,x) &\leq& - \left[ (1-\theta) D^2_{xx}c(0,y_0) + \theta D^2_{xx}c(0,y_1)\right]\cdot (x,x)\\
&& + \theta(1-\theta)|b-a|^2(- \demi C_0 |x'|^2 + \Delta_1|x_1||x|)\\
&\leq& - \left[ (1-\theta) D^2_{xx}c(0,y_0) + \theta D^2_{xx}c(0,y_1)\right]\cdot (x,x)\\
&& + \theta(1-\theta)|b-a|^2 (-\demi C_0 |x|^2 + (\Delta_1 + C_0)|x_1||x|).
\en
We set $\Delta_0 = \Delta_1 + C_0$.
Using a standard argument we have  $$\Delta_0|x||x^1| \leq C_0|x|^2/4 + |x^1|^2\Delta_0^2/C_0,$$ and we obtain
\be
-C_0 |x|^2/2 + \Delta_0|x||x_1|  \leq -C_0 |x|^2/4 + (\Delta_0^2/C_0) |x^1|^2.
\en
This concludes the proof of Lemma \ref{deltaDelta}.

$\hfill\Box$

\noindent
\textbf{Proof of Proposition \ref{long}.} Using the general fact that for  $f_0, f_1 \in \R$, for $0\leq \theta\leq 1$,
$$\max\{f_0,f_1\} \geq \theta f_1 + (1-\theta)f_0,$$  we have, using (\ref{ft}),
\be
\bar\phi (x) \geq &&(\theta b +(1-\theta)a)x^1\\
&-& \demi\left[ \theta D^2_{xx}c(0,y_1) + (1-\theta) D^2_{xx}c(0,y_0)\right]\cdot (x, x) \\
&+& o(|x|^2).
\en
We now use assumption {\bf As} through Lemma \ref{deltaDelta} to handle the second line of the right hand side. This yields the intermediate inequality
\beq\label{dodo}
\bar \phi(x) \geq && (\theta b +(1-\theta)a)x_1 -\demi D^2_{xx}c(0,y_\theta)\cdot (x,x) + \delta |x|^2  \\
&-& \Delta |x^1|^2\nonumber\\
&+& o(|x|^2),\nonumber
\enq
with $\delta, \Delta$ given in Lemma \ref{deltaDelta}.
In order to eliminate the term $-\Delta|x^1|^2$ in the right hand side, we proceed as follows:
We write  first (\ref{dodo}) for some $\theta' \in [0,1]$, and then  change it  into 
\be
\bar \phi(x) \geq && (\theta b +(1-\theta)a)x_1 -\demi D^2_{xx}c(0,y_\theta)\cdot (x,x) + \delta |x|^2 \\
&+& \demi\left[ D^2_{xx}c(0,y_\theta)-D^2_{xx}c(0,y_{\theta'})\right] \cdot (x,x)\\
&+&   ((b-a)(\theta'-\theta)-\Delta x^1)x^1 \\
&+&  (\delta'-\delta)|x|^2 + (\Delta - \Delta')|x^1|^2\\
&+& o(|x|^2), 
\en
where $\delta' = \delta(\theta'), \Delta'=\Delta(\theta')$ as in  Lemma \ref{deltaDelta}.
We now have to control the  terms
\be
T_1& = &((b-a)(\theta'-\theta)-\Delta x^1)x^1,\\
T_2&=& \demi\left[ D^2_{xx}c(0,y_\theta)-D^2_{xx}c(0,y_{\theta'})\right] \cdot (x,x),\\
T_3&=& (\delta'-\delta)|x|^2 + (\Delta - \Delta')|x^1|^2.
\en
The  term $T_1$ can be cancelled through an appropriate choice of $\theta'$.
We choose first  $\epsilon >0$ small (but fixed).
Taking $\theta \in [\epsilon, 1-\epsilon]$, we choose $\theta'$ such that 
\beq\label{defthetaprime}
\theta'= \theta+x^1 \Delta /(b-a)=\theta + x^1/C,
\enq
with 
\beq\label{defC}
C &=& (b-a)\Delta^{-1}\\
&=&  C_0 \left(\theta(1-\theta)|b-a| \Delta_0^2\right)^{-1}.\nonumber
\enq
Since $\theta \in [\epsilon, 1-\epsilon]$, this choice of $\theta'$ is possible if we restrict to $|x^1|\leq C\epsilon$. 
Note that $C_0$ is bounded away from 0, hence $C$ is bounded away from 0 for $|b-a|$ bounded (we don't need $C$ to be bounded from above).  Note also that $\theta'$ will depend on $x^1$ but $\theta$ has been fixed before. 

The second term $T_2$ is controlled using Lemma \ref{lemmeC1}: We have
\be
|T_2| \leq C_1|\theta'-\theta||x|^2.
\en

For the third term $T_3$, we note, from the definition of $\delta, \Delta$ in Lemma \ref{deltaDelta}, that $|\Delta - \Delta'|\leq |b-a|^2 \frac{\Delta_0^2}{C_0}|\theta'-\theta|$ and 
$|\delta - \delta'|\leq |C_0||\theta'-\theta|$.
Hence, using (\ref{defthetaprime}),
\be
|T_2 + T_3| 
&\leq& C_2|x|^3,
\en
where  $C_2$ depends on the bounds in  assumptions {\bf A0}, {\bf A2}, {\bf As}, and on $|b-a|$. We conclude that, for a suitable choice of $\theta'$,  
\beq
|T_1 + T_2+ T_3|\leq C_2|x|^3\label{estiT1T2T3}.
\enq 
We now have,  for all $\theta \in [\epsilon, 1-\epsilon]$, for all $x\in \Omega$ with $|x|<C\epsilon$, 
\be
\bar \phi(x) &\geq& (\theta b +(1-\theta)a)x_1 -\demi D^2_{xx}c(0,y_{\theta})\cdot (x, x)\\
&& + \delta |x|^2 \\
&& -  C_2|x|^3 +o(|x|^2).
\en
Using (\ref{ft}), this leads to
\be
\bar \phi(x) &\geq& -c(x, y_\theta) + c(0,y_\theta ) \\
&&+ \delta |x|^2 \\
&&- C_2|x|^3 +o(|x|^2).
\en 
We now notice that all the terms in $o(|x|^2)$ are error terms 
in the second order Taylor expansions (\ref{ft}).
Under assumption {\bf A0}, $c$ belongs to  $C^3(\bar\Omega\times\bar\Omega')$, hence there exists $\gamma$ such that the above inequality still holds true when replacing the third line of its right hand side by $-\gamma|x|^3$. The constant $\gamma$ will depend on the bounds in {\bf A0}, {\bf A2}, {\bf As}, and on $|b-a|$.
From Lemma \ref{deltaDelta}, we have $\delta = \frac{1}{4} C_0\theta(1-\theta)|b-a|^2$.
Using now (\ref{bma}), we have $$\frac{1}{C_\mT}|y_1-y_0| \leq |b-a|,$$ and letting $\delta_0=C_0 /(4C_\mT^2)$, 
we conclude the proof of Proposition \ref{long}.

$\hfill \Box$

\subsection{Part II. Construction of supporting functions}

We let ${\cal N}_{\mu}(B)$ denote the $\mu$-neighborhood of a set $B$, and we use the Proposition \ref{long} to prove the following:
\begin{prop}\label{supporting}
Let $\phi$ be c-convex. Let $c, \Omega, \Omega'$ satisfy the assumptions of Theorem \ref{main}. Let $x_0,x_1 \in \Omega$, and $y_0 \in G_\phi(x_0), y_1 \in G_\phi(x_1)$.
There exist constants $C, C', C''>0$ and $x_m\in [x_0,x_1]$, such that,
if
${\cal N}_\eta([x_0,x_1])\subset \Omega$, and 
\beq\label{pascontprop}
|y_1 - y_0| \geq \max\{|x_1-x_0|, C |x_1-x_0|^{1/5}\}>0,
\enq 
then 
$${\cal N}_\mu(\{y_\theta, \theta \in [1/4, 3/4]\})\cap \Omega' \subset G_\phi(B_\eta(x_m)),$$
where
\beq
&&\eta = C'\left(\frac{|x_1-x_0|}{|y_1-y_0|}\right)^{1/2},\label{defeta}\\
&&\mu= C'' \eta|y_1-y_0|^2\label{defmu}.
\enq
Here  $\{y_\theta\}_{\theta\in [0,1]}=[y_0,y_1]_{x_m}$ denotes the c-segment from $y_0$ to $y_1$ with respect to $x_m$.  Under assumptions {\bf A0}-{\bf As}, the constants $C, C'$ are bounded away from  infinity and $C''$ is bounded away from 0 .
\end{prop}

\textsc{Remark.} If  $x_0, x_1, y_0, y_1$ satisfying  (\ref{pascontprop}) can not be found, then $\phi$ is H\"older continuous with exponent $1/5$.

\paragraph{Preliminary result} Without loss of generality, we will assume that $\phi(x_0)=\phi(x_1)$: indeed, as remarked in Proposition \ref{subtract},
by subtracting from the cost function $c$ an affine function $\lambda$ that depends only on $x$, we will not modify the map solution of the optimal transportation problem, and the optimal potential $\phi$ will be changed into $\phi+\lambda$. Hence one can subtract a suitable affine function from $c$ so that $\phi(x_0)=\phi(x_1)$.
Notice that, as $\lambda$ is chosen affine, the gradient of the "new" potentials are deduced from the "old" ones just by adding the constant vector $\nabla_x \lambda$. Hence this does not change all the continuity properties of $\nabla \phi$, neither does it change all the derivatives of $c$ of order greater than or equal to 2.

As $y_0\in G_\phi(x_0), y_1 \in G_\phi(x_1)$ we have using  (\ref{def-c-transform}), for all $x \in \Omega$,
\be
-c(x,y_0) + c(x_0,y_0)+\phi(x_0) &\leq& \phi(x),\\
 -c(x,y_1) + c(x_1,y_1)+\phi(x_1) &\leq & \phi(x),
\en
with equality at $x=x_0$ in the first line, at $x=x_1$ in the second line.
Since $\phi(x_0)=\phi(x_1)$,
the difference between the supporting functions $x\to -c(x,y_0) + c(x_0,y_0)+\phi(x_0)$ and $x\to -c(x,y_1) + c(x_1,y_1)+\phi(x_1)$ will vanish  at some point $x_m$ in the segment $[x_0,x_1]$. Without loss of generality, we can add a constant to $\phi$ so that at this point both supporting functions are equal to 0.
Hence
\beq
&&-c(x_m,y_0) + c(x_0,y_0)+\phi(x_0)=0,\label{xm1}\\
&&-c(x_m,y_1) + c(x_1,y_1)+\phi(x_1)=0.\label{xm2}
\enq

\begin{lemme}\label{H}
Under the assumptions made above, and assuming moreover that $$|y_1-y_0| \geq |x_1-x_0|,$$ we have, for all $x$ in the segment $[x_0, x_1]$, 
$$\phi(x) \leq C_3|x_1-x_0||y_1-y_0|,$$
where $C_3$ depends only on $\|c(\cdot,\cdot)\|_{C^2(\Omega\times\Omega')}$.
\end{lemme}

\textsc{Proof of Lemma \ref{H}.} Using (\ref{xm1}, \ref{xm2}), we have 
\beq
H&=&\phi(x_0) \leq -\nabla_x c (x_m,y_0)\cdot (x_0-x_m) + \|c\|_{C^2}|x_0-x_m|^2/2\label{H2},\\
H&=&\phi(x_1) \leq  -\nabla_x c (x_m,y_1)\cdot(x_1-x_m) + \|c\|_{C^2}|x_1-x_m|^2/2\label{H3}.
\enq
By Proposition \ref{prop-c-convex}, the potential $\phi$ is semi-convex, with $D^2\phi \geq	-\|D^2_{xx}c\|_{\Linf(\Omega\times\Omega')}I$. Applying  the first part of Lemma \ref{semi-convex} to the function $f:t \to \phi(x_0 + t (x_1 - x_0))$ on $[0,1]$, for which $f''\geq -D^2\phi\cdot(x_1-x_0, x_1-x_0)$,  we find that  
\beq\label{H1}
\forall x\in[x_0,x_1], \;\phi(x) \leq H+C|x_1-x_0|^2,
\enq 
where $C = C(\|c\|_{C^2(\Omega\times\Omega')})$.
Then we consider two cases: 

The first one is where $-\nabla_x c (x_m,y_0)\cdot (x_0-x_m)$  and $-\nabla_x c (x_m,y_1)\cdot(x_1-x_m)$ are not both positive: let us assume for example that $-\nabla_x c (x_m,y_0)\cdot (x_0-x_m)$ is negative. Then we have, using (\ref{H2}), 
$H \leq \|c\|_{C^2}|x_0-x_m|^2/2$, and using (\ref{H1}), we get that
\be
\forall x\in[x_0,x_1], \;\phi(x) \leq (C + \|c\|_{C^2(\Omega\times\Omega')}/2)|x_1-x_0|^2.
\en 
Then we can conclude using $|x_1-x_0| \leq |y_1-y_0|$. 

We now consider the second case  where $-\nabla_x c (x_m,y_0)\cdot (x_0-x_m)$  and $-\nabla_x c (x_m,y_1)\cdot(x_1-x_m)$ are both positive. This implies that
\be
&& -\nabla_x c (x_m,y_0)\cdot (x_0-x_m) \leq -\nabla_x c (x_m,y_0)\cdot (x_0-x_1),\\
&& -\nabla_x c (x_m,y_1)\cdot (x_1-x_m) \leq -\nabla_x c (x_m,y_1)\cdot (x_1-x_0).
\en   
Combining with (\ref{H2}, \ref{H3}) we have
\be
2H &\leq&  -\nabla_x c (x_m,y_0)\cdot (x_0-x_1)  -\nabla_x c (x_m,y_1)\cdot(x_1-x_0) +\|c\|_{C^2}|x_0-x_1|^2 \\
&\leq& |\nabla_xc(x_m,y_0)-\nabla_xc(x_m,y_1)||x_0-x_1|+\|c\|_{C^2}|x_0-x_1|^2 \\
&\leq& \|c\|_{C^2}\Big(|x_1-x_0||y_1-y_0| +|x_0-x_1|^2\Big).
\en
Using  $|x_1-x_0| \leq |y_1-y_0|$, and then (\ref{H1}) we conclude.

$\hfill \Box$

\bigskip
We now assume the following: letting $\Gamma$ be defined by
\beq\label{Gamma}
\Gamma= \left[ \frac{\gamma^2}{\delta_0^3} 2^{12} C_3 \right]^{1/5}, 
\enq  
with $C_3(\|c\|_{C^2(\Omega\times\Omega')})$ defined in Lemma \ref{H},
$x_0, x_1, y_0, y_1$ satisfy
\beq\label{pascont}
|y_1-y_0| \geq \max\{\Gamma|x_1 - x_0|^{1/5}, |x_1 - x_0|\}.
\enq
Hence  the constant $C$ in Proposition \ref{supporting} will be equal to $\Gamma$.

\bigskip
We next state the following result, from which the proof of Proposition \ref{supporting} will follow easily:
\begin{lemme}\label{S}
 Let $x_m$ be defined as above.
For $y\in \Omega'$, consider the function  
\be
f_y(x)=-c(x,y)+c(x_m,y) + \phi(x_m).
\en
Under the assumptions made above, 
there exist $\eta, \mu$   as in Proposition \ref{supporting}, such that
for all $y\in {\cal N}_\mu(\{y_\theta, \theta \in [1/4, 3/4]\})\cap \Omega'$, 
\beq\label{eqS}
\phi - f_y \geq 0 \text{ on } \partial B_\eta(x_m)\cap\Omega.
\enq
\end{lemme}
Before proving this Lemma, we first show how it leads to Proposition \ref{supporting}.

\bigskip
\noindent
\textbf{Proof of Proposition \ref{supporting}.} By construction $f_y(x_m)=\phi(x_m)$, hence, if we have $\phi \geq f_y $ on $\partial B_\eta(x_m)$, then $\phi - f_y$ will have a local minimum inside $B_\eta(x_m)$, and for some point $x\in B_\eta(x_m)$, we will have $-\nabla_x c(x,y) \in \partial\phi(x)$. Using Theorem \ref{main2}, we have $\partial\phi(x) = \partial^c\phi(x)$, and this implies  $y\in G_\phi(x) \subset G_\phi(B_\eta(x_m))$. 

$\hfill\Box$

We now prove the main lemma:

\paragraph{Proof of Lemma \ref{S}.}
Using (\ref{xm1}, \ref{xm2}) and then  Proposition  \ref{long} centered at $x_m$ we obtain
\beq
\phi(x) 
&\geq& \max\{-c(x,y_0) + c(x_m,y_0),-c(x,y_1) + c(x_m,y_1)\}\nonumber\\
&\geq& -c(x,y_\theta) + c(x_m,y_\theta) + \delta_0 \theta(1-\theta)|y_0-y_1|^2|x-x_m|^2 -\gamma|x-x_m|^3\nonumber\\
&=&\Phi(x)\label{longxm}
\enq 
for $\epsilon >0$, for all $\theta\in [\epsilon,1-\epsilon]$, $|x-x_m| \leq C\epsilon$, and with $\{y_\theta\}_{\theta\in [0,1]}$ the c-segment with respect to  $x_m$ joining $y_0$ to $y_1$. 
Then we have for $y\in \Omega'$
\be
&&-c(x,y)+c(x_m,y) \\ 
&=& -c(x,y_\theta) + c(x_m,y_\theta) \\
&&+ \int_{s=0}^1 \left[\nabla_y c(x_m,y_\theta + s(y-y_\theta)) -\nabla_y c(x,y_\theta + s(y-y_\theta)) \right]  \cdot (y - y_\theta)\,ds \\
&\leq& -c(x,y_\theta) + c(x_m,y_\theta) + C_4|y-y_\theta||x-x_m|,
\en
where $C_4=\|D^2_{xy} c\|_{\Linf(\Omega\times\Omega')}$.
Combining this with Lemma \ref{H} to estimate $\phi(x_m)$, we have
\beq\label{fy}
f_y(x)&=&-c(x,y)+c(x_m,y) + \phi(x_m)\nonumber \\
&\leq& -c(x,y_\theta) + c(x_m,y_\theta) + C_4|y-y_\theta||x-x_m| + C_3|x_1-x_0||y_1-y_0|\nonumber\\
&=&F_y(x).
\enq
Inequality (\ref{eqS}) will be satisfied if we have, for $F_y, \Phi$ defined in (\ref{longxm}, \ref{fy})
\beq\label{eqSprime}
F_y(x) \leq \Phi(x)
\enq
on the set $\{|x-x_m|=\eta\}$, for some $\eta>0$.
First we restrict $\theta$ to $[1/4,3/4]$, (i.e. we take $\epsilon = 1/4$ in (\ref{longxm})). Then (\ref{eqSprime}) reads
\beq\label{whatwewant}
\frac{3}{16}\delta_0|y_0-y_1|^2 \eta^2  -\gamma  \eta^3  \geq C_4|y-y_\theta|\eta + C_3|x_1-x_0||y_1-y_0|.
\enq
Inequality (\ref{whatwewant}) will be satisfied if the three following inequalities are satisfied:
\be
\frac{1}{16}\delta_0|y_0-y_1|^2 \eta^2 & \geq &C_3|x_1-x_0||y_1-y_0|,\\
\frac{1}{16}\delta_0|y_0-y_1|^2 \eta^2 & \geq & C_4|y-y_\theta|\eta,\\
\frac{1}{16}\delta_0|y_0-y_1|^2 \eta^2 & \geq & \gamma  \eta^3.
\en
In order to satisfy the first inequality, we define $\eta$ by
\be
\eta^2 &=& \frac{16 C_3}{\delta_0} \frac{|x_1-x_0|}{|y_1-y_0|}.
\en
In order to satisfy the second, we define $\mu$ by
\be
\mu = C_5 \eta|y_1-y_0|^2,
\en
where $C_5 = \delta_0/(16C_4)$ (note that  $C_5$ is bounded away from 0), and consider $y\in \Omega'$ such that $|y-y_\theta| \leq \mu$.
The third inequality will then be implied by 
\be 
\gamma \eta \leq (\delta_0/16)|y_0-y_1|^2,
\en
which is equivalent to 
\be
\frac{\gamma^2}{\delta_0^3} 16^3 C_3 |x_1-x_0| \leq |y_1-y_0|^5, 
\en
and we recognize here assumption (\ref{pascont}). 
The constants $C, C', C''$ in Proposition \ref{supporting} are defined by  $C=\Gamma$ from assumption (\ref{pascont}), $C' = (\frac{16C_3}{\delta_0})^{1/2}$, $C'' = C_5$. Then, for all $y\in {\cal N}_\mu\{y_\theta, \theta\in [1/4,3/4]\}\cap \Omega'$, the function $f_y(x) = -c(x,y)+c(x_m,y) + \phi(x_m)$
will satisfy $f_y \leq \phi$ on  the boundary of the ball $B_\eta(x_m)$. 
This proves  Lemma \ref{S}.

$\hfill \Box$

\subsection{Part III. Continuity estimates}
 \begin{prop}\label{module} Let $\phi$ be c-convex with $G_\phi(\Omega)\subset \Omega'$. Let $c, \Omega, \Omega'$ satisfy the assumptions of Theorem \ref{main}.  Then, 
\begin{itemize}
\item 
if $G_\phi^{\#}{\rm dVol}$, satisfies (\ref{Lp}), for some $p>n$, then $\phi \in C^{1, \beta}_{loc}(\Omega)$, with $\beta(n,p)$ as in Theorem \ref{main},
\item
if $G_\phi^{\#}{\rm dVol}$ satisfies (\ref{mini}), then $\phi \in C^1_{loc}(\Omega)$,
\end{itemize}
where $G_\phi^\#$ is defined in Definition \ref{defiweak}.

\end{prop}

\paragraph{Preliminary Result}
We first state the following general result, whose proof is deferred to the appendix. 
\begin{lemme}\label{NinterOmega}
Let $\Omega'$ be c-convex with respect to $x_m\in \Omega$, let $y_0, y_1 \in \Omega'$. There exists $C, \mu_0 >0$ depending on $c, \Omega,  \Omega'$ such that for all $\mu \in (0, \mu_0)$, 
\be
{\rm Vol}\left( {\cal N}_\mu\left([y_0, y_1]_{x_m}\right)\cap\Omega' \right) \geq C {\rm Vol}\left( {\cal N}_\mu\left([y_0, y_1]_{x_m}\right)\right).
\en
\end{lemme}

\paragraph{Proof of Proposition \ref{module}}
Consider $\Omega_\delta = \{x\in \Omega, d(x,\partial\Omega) > \delta \}$.  In order to have ${\cal N}_\eta([x_0,x_1]) \subset \Omega$, it is enough to have
\begin{enumerate}
\item $x_0, x_1 \in \Omega_\delta$,
\item $|x_0-x_1| < \delta/2$,
\item $\eta < \delta /2$.
\end{enumerate}
If $y_i \in G_\phi(x_i), i=0,1$ satisfy (\ref{pascontprop}) in Proposition \ref{supporting}, $|y_1-y_0|\geq C|x_1-x_0|^{1/5}$,
then $\eta \leq E|x_1-x_0|^{2/5}$, with $\eta$ defined in  Proposition \ref{supporting}, and $E$ a constant depending only on $C', C''$ in Proposition \ref{supporting}.
 Hence for $|x_1-x_0|^{2/5} \leq \delta/(2E)$, it follows that ${\cal N}_\eta([x_0,x_1]) \subset \Omega$, and Proposition \ref{supporting} applies.
We now set 
\beq\label{defRdelta}
R_\delta = \inf\{\delta/2, (\delta/(2E))^{5/2}\}, 
\enq 
and in the remainder of the proof, we chose
$x_1, x_0 \in \Omega_\delta$  such that $|x_1-x_0| \leq R_\delta$.
From Proposition \ref{supporting}, we will have 
\beq\label{usesupp}
N_\mu\{y_\theta, \theta\in [1/4,3/4]\} \cap\Omega' \subset G_\phi(B_\eta(x_m)).
\enq
From Lemma \ref{NinterOmega}, and the definition of $\mu$ in (\ref{defmu}), there exits $C, C'>0$ such that 
\beq
{\rm Vol}\left(N_\mu\{y_\theta, \theta\in [1/4,3/4]\}\cap\Omega'\right) &\geq& C |y_1-y_0| \mu^{n-1}\label{volNmu}\\
&=& C'|y_1-y_0| \eta^{n-1}|y_1-y_0|^{2(n-1)}\nonumber.
\enq  

\paragraph*{$C^{1,\beta}$ estimates  for data with bounded density}
If the Jacobian determinant of the mapping $G_\phi$ is bounded, (in other words, if $G_\phi^\#{\rm dVol}$ has a density bounded in $\Linf$ with respect to  the Lebesgue measure) 
then, for some $C, C'$,   
\beq\label{volGphi}
{\rm Vol}\left(G_\phi(B_\eta(x_m))\right) &\leq& C{\rm Vol}\left(B_\eta(x_m)\right)\\
&=&C'\eta^n.\nonumber
\enq
Using (\ref{usesupp}) with (\ref{volNmu}), (\ref{volGphi}), we find for some $C,C'$ that
\be
|y_1-y_0|^{2n-1}&\leq& C \eta \\
&=& C'\left(\frac{|x_1-x_0|}{|y_1-y_0|}\right)^{1/2},
\en
which yields finally, for another constant $C_6>0$,
\be
|y_1-y_0|\leq C_6|x_1-x_0|^{\frac{1}{4n-1}}.
\en
From this we readily deduce  that $G_\phi$ is single valued, moreover 
$G_\phi \in C^{\frac{1}{4n-1}}_{loc}(\Omega)$. Since $-\nabla_xc(x,y_i) = \nabla\phi(x_i), i=0,1$, and $\nabla_xc$ is Lipschitz, this yields also $\phi \in C^{1, \frac{1}{4n-1}}_{loc}(\Omega)$.

\paragraph*{$C^{1,\beta}$ estimates for data satisfying (\ref{Lp})}
We can refine the argument: Let again $\nu = G_\phi^\#{\rm dVol}$ and $F$ be defined by 
\beq\label{defF}
F(V) &=& \sup\Big\{ {\rm Vol}(G_\phi(B)), B \subset \Omega \text{ a ball of volume } V\Big\}\\
&=& \sup\Big\{ \nu(B), B \subset \Omega \text{ a ball of volume } V\Big\}.\nonumber
\enq
Then, by Proposition \ref{supporting}, we have 
$F\left({\rm Vol}\left(B_\eta(x_m)\right)\right) \geq {\rm Vol}(N_\mu\{y_\theta, \theta\in [1/4,3/4]\}\cap\Omega')$, which yields, using (\ref{volNmu}) and the definition of $\eta$ in (\ref{defeta})
\beq\label{F}
F\left(\omega_n C_5^n\frac{|x_1-x_0|^{n/2}}{|y_1-y_0|^{n/2}}\right) \geq C_7|x_1-x_0|^{(n-1)/2}|y_1-y_0|^{(3n-1)/2}
\enq
for some  $C_7$  bounded away from 0, with $\omega_n$ the volume of the n-dimensional unit ball.
Assume that $F(V)\leq CV^\kappa$ for some $\kappa\in \R$. Note that $\nu\in L^p$ implies the (stronger) bound 
$F(V)=o(V^{1-1/p})$, hence it is natural to write $\kappa = 1-1/p$ for some $p\in ]1, +\infty]$, and the condition
\beq\label{howisF}
F(V) \leq C V^{1-1/p}
\enq
is then equivalent to condition (\ref{Lp}) for $\nu$.
We obtain from (\ref{F}) and (\ref{howisF}) that
\be
|y_1-y_0|^{2n-1 +\frac{1}{2}(1-\frac{n}{p})} \leq C_8 |x_1-x_0|^{\frac{1}{2}(1 - \frac{n}{p} )}.
\en
We see first that we need $p>n$, and, setting $\alpha = 1-n/p$, we obtain 
\be
|y_1-y_0| \leq C_9 |x_1-x_0|^{\frac{\alpha}{4n-2+\alpha}}.
\en
This yields H\"older continuity for $G_\phi$. Then we use that $\nabla\phi(x)=-\nabla_xc(x,G_\phi(x))$ and the smoothness of $c$ to obtain a similar H\"older estimate for $\nabla\phi$.

\paragraph{$C^1$ estimates for data satisfying (\ref{mini})}
We only  assume condition (\ref{mini}) for $\nu=G_{\phi}^{\#}{\rm dVol}$, which we can rewrite under the following form:
\beq\label{specialF}
F(V) \leq \left[f(V^{2/n})\right]^{2n-1}V^{1-1/n},
\enq for some  increasing $f:[0,1]\to \R^+$, with $\lim_{V\to 0} f(V) = 0$, $F$ being defined in (\ref{defF}).
Consistently with (\ref{pascont}), we can assume that, as $x_1$ goes to $x_0$, $\frac{|x_1-x_0|}{|y_1-y_0|}$
 goes also to 0. 
Using (\ref{specialF}) in (\ref{F}), we get for some $C_{10}, C_{11}$   bounded away from 0 and infinity,
\be
f^{2n-1}\left(C_{10}\frac{|x_1-x_0|}{|y_1-y_0|}\right) \geq (C_{11} |y_1-y_0|)^{2n-1},
\en
hence we get that $|y_1-y_0|$ goes to 0 when $|x_1-x_0|$ goes to 0. Then, let $g$ be the modulus of continuity of $G_\phi$ in $\Omega_\delta$; $g$ satisfies:

$\forall u \leq R_\delta$, either  $g(u) \leq \max\{u, \Gamma u^{1/5}\}$ or
\be  
&&f\left(C_{10}\frac{u}{g(u)}\right) \geq  C_{11}g(u),
\en
which is equivalent to 
\be
u \geq f^{-1}(C_{11}g(u))\frac{g(u)}{C_{10}},
\en
which is in turn equivalent to
\be
g(u) \leq \omega(u),
\en
where $\omega$ is the inverse of 
$\ds z\to f^{-1}(C_{11}z)\frac{z}{C_{10}}$. It is easily checked that $\lim_{r\to 0^+}\omega(r)=0$. This shows the continuity of $G_\phi$. 
Finally we have $\nabla\phi(x) = -\nabla_x c(x,G_\phi(x))$, and the continuity of $\nabla\phi$ is asserted.

$\hfill \Box$

\textsc{Remark.} The power $\beta=\frac{\alpha}{4n-2+\alpha}$ is not optimal for example if $n=1, p=+\infty$, for which the $C^{1,1}$ regularity is trivial, but note that in order to obtain this bound, we had to assume (\ref{pascont}). Hence the conclusion should be: either $\phi$ is $C^{1,1}$, or $\phi$ is $C^{1,1/5}$ or $\phi$ is $C^{1,\beta}$. Note that $\beta\leq 1/7$ for $n\geq 2$. 

\paragraph{Proof of Theorem \ref{main}}
In Proposition \ref{module}, we use a bound on $G_\phi^{\#}{\rm dVol}$. However, 
in Theorem \ref{main},  we only have $G_{\phi\,\#}\mu_1=\mu_0$, and as we we do not want to assume that $\mu_1\in L^1(\Rd)$, this does not imply necessarily that $G_\phi^{\#}\mu_1=\mu_0$ (See Definition \ref{defiweak} and the subsequent discussion). Hence we need the following proposition to finish the proof:
\begin{prop}\label{##}
Let $\phi$ be c-convex on $\Omega$, with $G_\phi(\Omega)\subset \Omega'$. Assume that $G_{\phi\,\#}\mu_0=\mu_1$. Assume that $\mu_1\geq m{\rm dVol}$ on $\Omega'$. Then for all $\omega \subset \Omega$, we have
\be
\mu_0(\omega) \geq m {\rm Vol}(G_\phi(\omega)),\text{ and hence, } G_\phi^{\#}{\rm dVol} \leq \frac{1}{m} \mu_0.
\en
\end{prop}

\textsc{Proof.} In $\Omega'$ we consider $N=\{y\in \Omega', \exists x_1\neq x_2 \in \Omega, G_\phi(x_1)=G_\phi(x_2)=y\}$. Then $N=\{y\in \Omega', \phi^c \text{ is not differentiable at } y\}$. Hence ${\rm Vol}(N)=0$, and  ${\rm Vol} (G_\phi(\omega)\setminus N) = {\rm Vol}(G_\phi(\omega))$. Moreover,  on $G_\phi(\omega)\setminus N$, $G_\phi^{-1}$ is single valued. Then  $G_\phi^{-1}(G_\phi(\omega)\setminus N) \subset \omega$. Hence, 
\be
\mu_0(\omega) &\geq& \mu_0 (G_\phi^{-1}(G_\phi(\omega)\setminus N))\\
&=&\mu_1(G_\phi(\omega)\setminus N)\\
&\geq & m{\rm Vol}(G_\phi(\omega)\setminus N)\\
& = & m{\rm Vol}(G_\phi(\omega)).
\en 

$\hfill \Box$

\paragraph{Proof of the boundary regularity}
This part is easy: under the assumptions of Theorem \ref{main4}, the density $\mu_0$ satisfies (\ref{Lp}) with $p>n$ (resp. satisfies (\ref{mini})). Hence Theorem \ref{main} applies and $\phi \in C^{1,\beta}_{loc}(\Omega)$ (resp.  $\phi \in C^1_{loc}(\Omega)$). Since $\Omega_2$ is compactly contained in $\Omega$, we conclude the boundary regularity on $\Omega_2$. 
 This proves Theorem \ref{main4}. $\hfill \Box$

\textsc{Remark.} This proof of the boundary regularity is very simple because we have interior regularity even when $\mu_0$ vanishes. This is not the case for the classical Monge-Amp\`ere equation, and the boundary regularity requires that both $\Omega$ and $\Omega'$ are convex, and is more complicated to establish (see \cite{Ca4}).

We now show that there is indeed equivalence between assumption {\bf As} at a point $x$ and the conclusion of Proposition \ref{long}. This is a quantitative version of Theorem \ref{anoune}.
\begin{prop}\label{reciplong}
Assume that at a point $x_0$ for all $y_0, y_1$, for $y_{1/2}$ the 'middle' point of $[y_0,y_1]_{x_0}$, we have
$$\bar\phi(x) \geq -c(x,y_{1/2})+c(x_0,y_{1/2}) + \delta_0 |y_0-y_1|^2|x-x_0|^2 +O(|x-x_0|^3)$$
with $\bar\phi$ as above. Then the cost function satisfies assumption {\bf As} at $x_0$ with $C_0=C\delta_0$, for some constant $C>0$ that depends on the bound in {\bf A2}.
\end{prop}

\textsc{Proof.} The proof follows the same lines as the proof of Theorem \ref{main2}, and is omitted here.

$\hfill\Box$

\section{Proof of Theorem \ref{costcurvature}}

 We consider condition {\bf Aw} at $(x_0,y=x_0)$.  We recall that
\be
\mathfrak{S_c}(x_0,x_0)(\xi,\nu)=-D^2_{p_\nu p_\nu}D^2_{x_\xi x_\xi} [(x,p)\to c(x, \mathfrak{T}_{x_0}(p))].
\en
for any $\nu,\xi$ in $T_{x_0}M$.
Let us first take a normal system of coordinates at $x_0$, so that we will compute
$$
Q=-D^2_{tt}D^2_{ss} [(x,p)\to c(\mathfrak{T}_{x_0}(t\xi), \mathfrak{T}_{x_0}(s\nu))].
$$
Let us write a finite difference version of this operator.
We first introduce $y_- = \mathfrak{T}_{x_0}(-h \nu), y_+= \mathfrak{T}_{x_0}(h\nu), x_-=\mathfrak{T}_{x_0}(-h \xi), x_+=\mathfrak{T}_{x_0}(h \xi)$.
We use the usual second order difference quotient, for example
\be
D^2_{x_\xi,x_\xi}c(x,\mathfrak{T}_{x_0}(p)) = \lim_{h\to 0}\frac{1}{h^2}(c(x_+,x_0)-2c(x_0,x_0)+c(x_-,x_0)).
\en
(Of course we have $c(x_0,x_0)=0$.)
We will have, as $h$ goes to $0$,
\be
&&\lim_{h\to 0}\frac{1}{h^4}\left(\sum_{i,j=+,-}c(x_i, y_j) -2\sum_{i=+,-}(c(x_i,x_0)+c(y_j,x_0))\right)=-Q.\\
\en
 Rearranging the terms, we find that the left hand side of the previous identity is equal to 
\be
\sum_{i,j=+,-}[c(x_i, y_j)-c(x_0,x_i)-c(x_0,y_j)].
\en
Each of the terms inside brackets has a simple geometric interpretation: consider the triangle with vertices $(x_0, x_i, y_j)$ whose sides are geodesics. This is a square angle triangle. If the metric is flat, by Pythagoras Theorem, the term inside the brackets is 0. In the general case, a standard computation shows that it is equal to $-\frac{1}{6}\kappa(x_0,\xi, \nu) h^4 + o(h^4)$ where $\kappa(x_0, \xi,\nu)$ is the sectional curvature at $x_0$ in the two-plane generated by $\xi, \nu$. Hence, we get that $Q=(2/3)\kappa(x, \xi,\nu)$.

Now to reach the more general formula of Theorem \ref{costcurvature}, we use the following expansion of the distance that C\'edric Villani communicated us:
\begin{lemme}
Let $M$ be a smooth Riemannian manifold. Let $\gamma_1,\gamma_2$ be two unit speed geodesics that leave point $x_0\in M$. Let $\theta$ be the angle between $\dot\gamma_1(0)$ and $\dot \gamma_2(0)$ (measured with respect to the metric), let $\kappa$ be the sectional curvature of $M$ at $x_0$ in the 2-plane generated by $\dot\gamma_1(0), \dot\gamma_2(0)$.
Then we have
\be
d^2(\gamma_1(t), \gamma_2(t))=2(1-\cos(\theta))(1-\frac{\kappa}{6}(\cos^2(\theta/2)) t^2 + O(t^4))t^2 )^2.                   
\en
\end{lemme}
Then, we obtain easily, following the same lines as in the case looked above that
\be
\mathfrak{S_c}(x_0,x_0)(\xi,\nu)
=(2/3)\kappa(x_0, \xi,\nu) (|\xi|^2_g|\nu|^2_g - (\xi,\nu)_g^2),
\en
where $(\cdot,\cdot)_g, |\cdot|_g$ denote respectively the scalar product and the norm with respect to $g$.
This proves the Theorem. $\hfill\Box$

\subsection{Counterexample to regularity for a manifold with negative curvature}
Consider the two dimensional surface $H=\{z=x^2-y^2\}\subset \R^3$, endowed with the Riemannian metric inherited from the canonical metric of $\R^3$. Then $H$ has negative sectional curvature around $0$. For $r$ sufficiently small, $\Omega=H\cap B_r(0)$ is c-convex with respect to itself. Consider the function 
\be
\bar\phi(x)=\max\{-d^2/2(X,X_0), -d^2/2(X,X_1)\},
\en
where $X_0=(0,a,-a^2), X_1=(0,-a, -a^2)$. Then, as shown by our proof of Theorem \ref{main2}, for $a$ small enough, no sequence of $C^1$ c-convex potentials can converge uniformly to $\bar\phi$ on $\Omega$.
Let $\mu_0$ to be the Lebesgue measure of $\Omega$, and $\mu_1=\demi(\delta_{X_0}+\delta_{X_1})$. We have $G_{\bar\phi\,\#}\mu_0 = \mu_1$. Let $\mu_1\ep \in C^\infty(\bar\Omega)$  be a positive mollification of $\mu_1$ so that its total mass remains equal to 1, and that preserves the symmetries with respect to $x=0$ and $y=0$. 
Let $\phi_n$ be such that $G_{\phi_n\,\#}\mu_0=\mu_n$.
Then, for $n$ large enough, $\phi_n$ is not differentiable at the origin. Indeed, for symmetry reasons, 0 belongs to the subdifferential of $\phi_n$ at 0, on the other hand,  $\phi_n$ converges uniformly to $\bar\phi$, and we know from the fact that {\bf Aw} is violated at 0 that $-\bar\phi-c(\cdot,0)$ does not reach its global maximum on $\Omega$ at 0.

\section{Appendix}
\paragraph{Proof of proposition \ref{propnewR}}
We first prove the ``intrinsic'' part. 
In order to show this, we consider $\gamma$ a $C^2$ curve in $\Omega$ defined in a neighborhood of 0, such that 
\beq
\label{1ga}\gamma(0) &=& x_0,\\ 
\label{2ga}\dot\gamma(0)&=& \xi.
\enq
We then consider the quantity
\be
Q_\gamma = D^2_{tt}D^2_{ss} \Big[(s,t)\to c(\gamma(t), \mathfrak{T}_{x_0}(p_0+s\nu))\Big]\Big|_{t,s=0}.
\en
We show that this quantity is independent of the choice of $\gamma$.
We have
\be
Q_\gamma &=& D^2_{ss}\Big(s\to D^2_{\xi\xi}c(x_0, \mathfrak{T}_{x_0}(p_0+s\nu))+ D_x c(x_0,\mathfrak{T}_{x_0}(p_0+s\nu))\cdot\ddot\gamma(0)\Big)\\
&=& D^2_{ss}\Big(s\to D^2_{\xi\xi}c(x_0, \mathfrak{T}_{x_0}(p_0+s\nu))- (p_0+s\nu)\cdot\ddot\gamma(0)\Big)\\
&=& D^2_{ss}\Big(s\to D^2_{\xi\xi}c(x_0, \mathfrak{T}_{x_0}(p_0+s\nu))\Big),
\en
where the second line follows from the very definition of the c-exponential map.
Hence, the value of the curvature is independent of $\ddot\gamma(0)$, and therefore of the choice of $\gamma$ as long as it satisfies (\ref{1ga}, \ref{2ga}). One can now choose around $x_0$ a system of geodesic coordinates, which yields the equivalence of the definitions (\ref{defR}) and (\ref{defRint}).
Then, the second part of Proposition \ref{propnewR} follows by taking as new coordinates around $x_0$ the c-geodesics with respect to $y_0$, which yields

\beq\label{newnewRint}
\mathfrak{S_c}(x_0,y_0)(\xi,\nu) = D^2_{p_\nu p_\nu q_{\tilde\xi} q_{\tilde\xi}}\left[c(\mathfrak{T}_{y_0}(q),\mathfrak{T}_{x_0}(p)) \right]\Big|_{q_0=-\nabla_y c(x_0,y_0),\,p_0=-\nabla_x c(x_0,y_0)},
\enq
where $\tilde \xi$ is chosen such that 
\be
D_{q}\mathfrak{T}^*_{y_0}(q_0)\cdot\tilde\xi=\xi.
\en
The condition $\xi\perp\nu$ nows reads $(D_{q}\mathfrak{T}^*_{y_0}(q_0)\cdot{\tilde\xi}) \perp\nu$ or equivalently
$[D_{x,y}c]^{-1}\cdot(\nu,\tilde\xi)=0$.
Then, identity (\ref{newRint}) follows by a symmetric argument.

$\hfill \Box$

\paragraph{Proof of Proposition \ref{proprho}.} We prove only the last point, the other points being elementary. Consider on $\Rd$ a measure locally equal to $\mu_0={\cal L}^{n-1} \otimes \mu$, where ${\cal L}^{n-1}$ is the $n-1$-dimensional Lebesgue measure, and $\mu$ is a probability measure on $[0,1]$ equal to the derivative of the Devil's staircase. Then, $\mu \notin L^1$. On the other hand, for all $[a,b]\subset[0,1]$, $\mu([a,b]) \leq |b-a|^\alpha$, for some $\alpha \in (0,1]$. Then, for $x=(x_1,..,x_n)$,  
$$\mu_0(B_r(x))\leq C r^{n-1}\mu[x_n-r,x_n+r]) \leq C r^{n-1+\alpha} = Cr^{n(1-1/p)}$$
for some $p>n$. Hence $\mu_0 \notin L^1_{loc}$ and $\mu_0$ satisfies (\ref{Lp}) for some $p>n$.

\paragraph{Proof of Proposition \ref{unicite}} We know (see \cite[chapter 2]{Vi}) that there exists $\pi$ a probability measure on  $\Rd\times\Rd$, with marginals $\mu$ and $\nu$, and such that
\be
\int_{\Rd} \phi(x) d\mu(x) + \psi(x)d\nu(x) = -\int c(x,y)d\pi(x,y),
\en and moreover, there exists $\bar\phi$ a c-convex potential such that $$supp(\pi)\subset\{(x,G_{\bar\phi} (x)), x\in \Rd\}.$$ Let us decompose $\pi$ as $\pi=\mu\otimes\gamma_x$, where for $d\mu$ almost all $x\in \Rd$ ,  $\gamma_x$ is a probability measure on $\Rd$ and $\gamma_x$ is supported in $G_{\bar\phi}(x)$.
Hence we have
\be
\int_{\Rd}d\mu(x)\left[\int_{\Rd}d\gamma_x(y)(\phi(x)+\psi(y)-c(x,y))\right] = 0.
\en
This implies that for $d\mu$ a.e. $x$, for $d\gamma_x$ a.e. $y$, we have $y \in G_{\phi}(x)$. Since for $d\mu$ a.e. $x$, we have $y\in G_{\bar\phi}(x)$ $d\gamma_x$ a.s., we deduce that for $d\mu$ a.e. $x$, (and hence for Lebesgue a.e. $x$, since $\mu>0$ a.e.), we have $G_{\bar\phi}(x)\cap G_{\phi}(x)\neq \emptyset$. This implies that
$\nabla\phi=\nabla\bar\phi$ Lebesgue a.e., and that $\phi-\bar\phi$ is constant. 
This shows that $\phi$ is uniquely defined up to a constant.
Now the pair $\psi^{c*}, \psi$ can only improve the infimum (\ref{kanto}) compared to $(\phi, \psi)$, hence it is also optimal. Hence $\psi^{c*}$ is also uniquely defined up to a constant. If $\psi$ is c*-convex, then $\psi^{c*c}=\psi$, and $\psi$ is thus uniquely defined.

$\hfill\Box$

\paragraph{Proof of Lemma \ref{NinterOmega}}
From {\bf A1}, {\bf A2}, for all $x_m\in \Omega$, $\psi: y\to -\nabla_xc(x_m,\cdot)$ is a diffeomorphism from $\Omega'$ to $-\nabla_xc(x_m,\Omega')$. Then 
\be
\psi \left( {\cal N}_\eta\left([y_0, y_1]_{x_m}\right)\cap\Omega' \right) = \psi\left( {\cal N}_\eta\left([y_0, y_1]_{x_m}\right) \right) \cap \psi\left(\Omega'\right).
\en
Letting $p_i = -\nabla_xc(x_m, y_i), i=0,1$, using {\bf A1}, {\bf A2}, there exists $C>0$ such that
\be
{\cal N}_{C\eta} \left([p_0, p_1]\right) \subset \psi\left( {\cal N}_\eta\left([y_0, y_1]_{x_m}\right) \right).
\en 
Moreover, as $\Omega'$ is c-convex with respect to $x_m$, $\psi(\Omega')$ is a convex set.

Then we claim the following: for $U\subset\Rd$ convex, for $u,v\in U$, the function
\be
r \to {\rm Vol}({\cal N}_r([u,v]) \cap U)/{\rm Vol}({\cal N}_r([u,v]))
\en is non-increasing. Indeed, by convexity of $U$, for $w \in [u,v]$, if $w+w' \in B_r(w) \cap U$, for $\theta \in [0,1]$, $w+ \theta w' \in B_{\theta r}(w)\cap U$. Then the claim follows easily.

Hence, we have 
\be
&&{\rm Vol}\left(\psi\left( {\cal N}_\eta\left([y_0, y_1]_{x_m}\right)\cap\Omega' \right)\right) \\
&\geq & {\rm Vol}\left( {\cal N}_{C\eta} \left([p_0, p_1]\right) \cap \psi(\Omega') \right)\\
&\geq & {\rm Vol}\left( {\cal N}_{C\eta} \left([p_0, p_1]\right)\right)  {\rm Vol}\left( {\cal N}_{1} \left([p_0, p_1]\right) \cap \psi(\Omega') \right) {\rm Vol}^{-1}\left( {\cal N}_{1} \left([p_0, p_1]\right)\right),
\en
whenever $\eta$ is small enough so that $C\eta \leq 1$.
By compactness, one has
\be
{\rm Vol}\left( {\cal N}_{1} \left([p_0, p_1]\right) \cap \psi(\Omega') \right) {\rm Vol}^{-1}\left( {\cal N}_{1} \left([p_0, p_1]\right)\right) \geq C(\Omega').
\en
Moreover, for $C>0$, there exists a constant $C'>0$ such that
\be
{\rm Vol}\left( {\cal N}_{C\eta} \left([p_0, p_1]\right)\right) \geq C' {\rm Vol}\left( {\cal N}_{\eta} \left([p_0, p_1]\right)\right)
\en
for all $\eta>0$.
Then, as $\psi$ is a smooth diffeomorphism, one has
\be
&&{\rm Vol}\left( {\cal N}_\eta\left([y_0, y_1]_{x_m}\right)\cap\Omega' \right) / {\rm Vol}\left( {\cal N}_\eta\left([y_0, y_1]_{x_m}\right)\right)\\
&\geq& C(c, \Omega, \Omega') {\rm Vol}\left( \psi\left({\cal N}_\eta\left([y_0, y_1]_{x_m}\right)\cap\Omega' \right)\right) / {\rm Vol}\left(\psi \left( {\cal N}_\eta\left([y_0, y_1]_{x_m}\right)\right)\right).
\en

$\hfill\Box$

\bibliography{C1-biblio}
\vspace{1cm}
\begin{flushright} 
Gr\'egoire Loeper
\\
Universit\'e Claude Bernard Lyon 1
\\ 
43 boulevard du 11 novembre 1918
\\
69622 Villeurbanne cedex France
\\
loeper@math.univ-lyon1.fr
\end{flushright}

\end{document}